\documentclass[12pt,reqno]{amsart}
\usepackage{mathrsfs}
\usepackage{cases}
\usepackage{epic}
\usepackage{amsfonts}
\usepackage{graphicx}
\usepackage{amsmath}
\usepackage{amssymb, upgreek}
\usepackage{bm}
\usepackage{latexsym,todonotes}
\usepackage{pdflscape}
\usepackage[all]{xypic}
\usepackage[all]{xy}
\usepackage{color}
\usepackage{colordvi}
\usepackage{multicol}

\usepackage[linktocpage=true]{hyperref}
\hypersetup{colorlinks,linkcolor=blue,urlcolor=cyan,citecolor=blue}
\usepackage{tikz}
\newcommand{\nc}{\newcommand}
\nc{\browntext}[1]{\textcolor{brown}{#1}}
\nc{\greentext}[1]{\textcolor{green}{#1}}
\nc{\redtext}[1]{\textcolor{red}{#1}}
\nc{\bluetext}[1]{\textcolor{blue}{#1}}
\nc{\brown}[1]{\browntext{ #1}}
\nc{\green}[1]{\greentext{ #1}}
\nc{\red}[1]{\redtext{ #1}}
\nc{\blue}[1]{\bluetext{ #1}}
\usepackage{latexsym,todonotes}
\usepackage{tipa}

\newtheorem{theorem}{Theorem}  [section]
\newtheorem{corollary}[theorem]{Corollary}
\newtheorem{lemma}[theorem]{Lemma}

\newtheorem{proposition}[theorem]{Proposition}

\theoremstyle{remark}
\newtheorem{remark}[theorem]{Remark}

\numberwithin{equation}{section}

\newcommand{\cc}{{\mathcal C}}

\def \co{\mathcal O}

\def \ch{{\mathcal H}}

\def \ct{{\mathcal T}}

\numberwithin{equation}{section}

\renewcommand{\ker}{\operatorname{Ker}\nolimits}

\def \cp{\mathcal P}
\newcommand{\Hom}{\operatorname{Hom}\nolimits}

\newcommand{\Aut}{\operatorname{Aut}\nolimits}

\newcommand{\Ext}{\operatorname{Ext}\nolimits}

\def \och{\check{\ch}}

\newcommand{\mbf}{\mathbf}

\newcommand{\mrm}{\mathrm}

\newcommand{\End}{\mrm{End}}
\newcommand{\rank}{\mrm{rank}}

\def \C{{\mathbb C}}

\newcommand{\la}{\lambda}

\newcommand{\Q}{\mathbb Q}

\newcommand{\U}{\mbf U}

\newcommand{\F}{\mathbb F}

\newcommand{\arxiv}[1]{\href{http://arxiv.org/abs/#1}{\tt arXiv:\nolinkurl{#1}}}

\newcommand{\Z}{\mathbb Z}

\def \X{\mathbb X}

\def \cv{\mathcal V}

\def \bM{\textbf{\Large\textturnw}}

\newcommand{\sqq}{{\bf v}}

\newcommand{\coker}{\operatorname{Coker}\nolimits}

\def \cc{\mathcal C}
\def \cf{\mathcal F}
\def \cg{\mathcal G}

\def \tt{\widetilde{t}} 

\def\ca{\mathcal A}
\newcommand{\Iso}{\operatorname{Iso}\nolimits}
\renewcommand{\Im}{\operatorname{Im}\nolimits}

\newcommand{\coh}{\operatorname{coh}\nolimits}
\newcommand{\rad}{\operatorname{rad}\nolimits}
\newcommand{\rep}{\operatorname{rep}\nolimits}
\newcommand{\Ker}{\operatorname{Ker}\nolimits}

\newcommand{\cone}{\operatorname{Cone}\nolimits}

\def \cw{\mathcal W}

\def \ch{\mathcal H}
\def \cd{\mathcal D}
\def\bfk{\mathbf{k}}

\def \QJ{Q_{\texttt{J}}}

\def \th{\widetilde{\ch}}

\def \cR{\mathcal R}
\def \ce{\mathcal E}
\def \cb{\mathcal B}

\def \bx{\mathbf{x}}

\allowdisplaybreaks

\begin{document}

	\title[Derived Hall algebras of root categories]
	{Derived Hall algebras of root categories}

	\author[Jiayi Chen]{Jiayi Chen}
	\address{School of Mathematical Sciences,
		Xiamen University, Xiamen 361005, P.R.China}
	\email{jiayichen.xmu@foxmail.com}
	
	\author[Ming Lu]{Ming Lu}
	\address{Department of Mathematics, Sichuan University, Chengdu 610064, P.R.China}
	\email{luming@scu.edu.cn}
	
	\author[Shiquan Ruan]{Shiquan Ruan}
	\address{ School of Mathematical Sciences,
		Xiamen University, Xiamen 361005, P.R.China}
	\email{sqruan@xmu.edu.cn}

	
	
	\begin{abstract}

		For a finitary hereditary abelian category $\ca$, we define a derived Hall algebra of its root category by counting the triangles and using the octahedral axiom, which is proved to be isomorphic to the Drinfeld double of Hall algebra of $\ca$. When applied to finite-dimensional nilpotent representations of the Jordan quiver or coherent sheaves over elliptic curves, these algebras provide categorical realizations of the ring of Laurent symmetric functions and also double affine Hecke algebras.
	\end{abstract}
	
	\maketitle
	\setcounter{tocdepth}{1}
	\tableofcontents

	
	\section{Introduction}
	
	\subsection{Backgrounds}

	Hall algebra is a fundamental concept in representation theory with intimate connections and applications to combinatorics of symmetric polynomials,  quantum groups and algebraic geometry.  The first occurrence of the concept of Hall algebras can
	be dated back to the early days of the twentieth century in the work of
	Steinitz (systematically studied by Hall in 1950s) which deals
	with the case of the category $\ca$  of abelian $p$-groups for $p$  a fixed prime number, or equivalently the category $\ca$ of finite-dimensional nilpotent representations of the Jordan quiver. 
	This classic Hall algebra is isomorphic to the ring $\Lambda_t$ of symmetric polynomials, and its Hall basis corresponds to Hall-Littlewood symmetric functions, a distinguished basis of $\Lambda_t$; cf. \cite{Mac95}.   
	
	Following Steinitz and Hall, Ringel \cite{Rin90} in 1990 defined Hall algebras for general (hereditary) abelian categories, which provided a realization of half parts of quantum groups and led to Lusztig's construction of the canonical basis \cite{Lus93}. Later on, Green \cite{Gr95} and Xiao \cite{X97} constructed the coproduct and the antipode for the Hall algebra of a hereditary abelian category $\ca$, and then Xiao used the Drinfeld double of Hall algebra to realize the whole quantum group. However, the Drinfeld double of Hall algebra is not so perfect, since it is built of two copies of (extended) Hall algebras  and the multiplication between the two pieces has to be put in by hand, so its construction is not purely categorical.

	Since then, many experts tried to realize the whole quantum groups by using Hall algebras.  In 2006, To\"{e}n \cite{T06} introduced the derived Hall algebra for a dg-category (under some finiteness conditions), generalizing the Hall algebra of an abelian
	category. Later, Xiao and Xu \cite{XX08} generalized To\"{e}n's construction to any triangulated categories with
	(left) homological-finite condition, which include the bounded derived category of a finitary hereditary abelian category $\ca$. However, as stated in many articles (such as \cite{T06,XX08,XC13}), the full quantum
	group is related not to the bounded derived category of $\ca$ but rather to its root category.  The root category of $\ca$ does not satisfy the left homological-finite condition, so its derived Hall algebra is not defined by their work. 
	The root category is a periodic triangulated category, which motivated Xu and Chen \cite{XC13} to construct Hall algebras for odd periodic triangulated categories. 
	
	Bridgeland \cite{Br13} in 2013 constructed a so-called Bridgeland's Hall algebra by taking localization of the Hall algebra of the category of 2-periodic projective complexes. He proved that  the whole quantum group can be realized by this new Hall algebra. %
	Bridgeland's Hall algebra has been found to have further generalizations and improvements which allow more flexibilities. Gorsky \cite{Gor18} constructed  semi-derived Hall algebras for Frobenius categories. Motivated by the works of Bridgeland and Gorsky, Lu and Peng \cite{LP21} formulated the  semi-derived (Ringel-)Hall algebras by using $2$-periodic complexes of hereditary abelian categories. Recently, Lu and Wang \cite{LW22,LW20} introduced $\imath$quiver algebras and used their semi-derived Hall algebras to realize the universal $\imath$quantum groups (also called quantum symmetric pair coideal subalgebras).  Based on \cite{LW22,LW20}, the authors  \cite{CLR22} used the derived Hall algebra of 1-periodic derived category \cite{XC13} to realize the original $\imath$quantum group of split type defined in \cite{Let99,Ko14}.

	By Bridgeland’s work \cite{Br13}, it seems
	that the exact category of 2-periodic complexes is the right setting rather than the root category to
	be related to the quantum group. However, as we shall see, the root categories have deep connections to  some other interesting rings such as the ring of Laurent symmetric functions and double affine Hecke algebra (DAHA). This motivates us to consider the open question of defining derived Hall algebras for root categories.
	
	The Hall algebra and also its Drinfeld double of an elliptic curve $\X$ over a finite field have been studied by many experts (see \cite{BS12,S12,SV1}), and were found to have deep relationships to Macdonald's polynomials and DAHAs defined in \cite{Ch05}. In \cite{BS12}, Burban and Schiffmann constructed a subalgebra $\U_\X^+$ of the Hall algebra of the category $\coh(\X)$ of coherent sheaves over $\X$ to give a categorical realization of the ring 	$\bM^+=\C[x_1^{\pm},x_2^{\pm},\dots,y_1,y_2,\dots]^{\mathfrak{S}_\infty}$. Then they used the Drinfeld double $\U_\X$ of $\U_\X^+$ to realize the ring of diagonal symmetric polynomials $\bM=\C[x_1^{\pm},x_2^{\pm},\dots,y_1^\pm,y_2^\pm,\dots]^{\mathfrak{S}_\infty}$. Besides, the ring $\Lambda_t$ of symmetric polynomials is a Hopf algebra, and its Drinfeld double $\cd\Lambda_t$ (also called the ring of Laurent symmetric functions) is also interesting and important in representation theory and combinatorics; see \cite{BS12,SV15,Zel81}. We can also use the Drinfeld double of Hall algebra of the Jordan quiver to realize $\cd\Lambda_t$ based on the works of Hall, Steinitz and Macdonald. 
	
	It is a natural question to give a categorical realization of  $\bM$ and $\cd\Lambda_t$ via Hall algebras. 
	However, the semi-derived Hall algebras of the elliptic curves and Jordan quiver can not be used {\em directly} to realize these interesting rings. For example, the semi-derived Hall algebra of the Jordan quiver is not commutative but $\cd\Lambda_t$ is. As we see, the root categories are the right setting to construct Hall algebras to realize these  interesting rings directly.

	\subsection{Main results}
	
	The goal of this paper is to formulate derived Hall algebras for root categories of hereditary abelian categories. 
	
	Let $\ca$ be a hereditary abelian category with skew symmetric Euler form. 
	By counting the triangles in root category $\cR(\ca)$ of $\ca$, we define the derived Hall number $F_{L^\bullet,Z^\bullet}^{M^\bullet}$ in \eqref{3.9} for any objects $L^\bullet,M^\bullet,Z^\bullet$, which also involves the Euler form of $\ca$. The derived Hall algebra  is a $\Q$-space with a basis parameterized by the isomorphism classes $[X^\bullet]$ of the root category, and the derived Hall numbers being its structure constants. We prove that the derived Hall algebra is an associative algebra by using the octahedral axiom, and then prove that it is isomorphic to the Drinfeld double of Hall algebra of $\ca$.  The category of finite-dimensional nilpotent representations of the Jordan quiver and the categories of coherent sheaves over elliptic curves are typical examples.
	As applications, we study the derived Hall algebras of the Jordan quiver and elliptic curves, which are used to realize the ring of Laurent symmetric functions and DAHAs. 

	However, one can also define derived Hall algebras for arbitrary hereditary abelian categories by using the same method of this paper. Note that the construction of derived Hall algebras in this general setting could  not  be direct since  a quantum torus (i.e., the group algebra of the Grothendieck group $K_0(\ca)$) should be appended; see Remark \ref{rem}. As a result, the derived Hall algebras constructed in this way are isomorphic to the semi-derived Hall algebras \cite{LP21} and also the Drinfeld double of extended Hall algebras (see \cite{X97}), and can be used to realize the quantum groups; cf. \cite{X97,LP21,LW22}. 

	\subsection{Comparison with previous works}
	We know that the precise formula of derived Hall numbers \eqref{3.9} (which equals to \eqref{3.10}--\eqref{3.11}) is crucial to define the derived Hall algebra.
	The derived Hall number defined in  \cite[\S3]{XX08} is mainly to compute the cardinality of the set $\Hom(M^\bullet,L^\bullet)_{Z^\bullet[1]}=\{l:M^\bullet\rightarrow L^\bullet\mid \cone(l)\cong Z^\bullet[1]\}$. In order to give the derived Hall numbers here, we need to divide the set 
	$\Hom(M^\bullet,L^\bullet)_{Z^\bullet[1]}$ into subsets  $\Hom(M^\bullet,L^\bullet)_{Z^\bullet[1],\delta}$ for $\delta\in K_0(\ca)$ (see \eqref{eq:HomZMLdelta}). The derived Hall number is to compute the cardinalities of these subsets and sum over $\delta\in K_0(\ca)$  by multiplying suitable Euler forms. 
	The idea to prove the associativity of the derived Hall algebra is almost the same as \cite{XX08} but with some much subtle  computations.
	
	Recently, Zhang \cite{Zh22} defined a Hall algebra
	for the root category of $\ca$ by applying the derived Hall numbers of the bounded derived
	category $\cd^b(\ca)$, which is proved to be isomorphic to the Drinfeld double of Hall algebra of $\ca$. Our derived Hall algebra is certainly isomorphic to Zhang's, but these two definitions are much different: the structure constants of Zhang's Hall algebra are to count the triangles in $\cd^b(\ca)$, and ours are to count the triangles in the root category of $\ca$.

	\subsection{Future works}
	
	For the Jordan quiver, as in \S\ref{sec:Jordan}, its derived Hall algebra could be used to realize the ring of Laurent symmetric functions $\cd\Lambda_t$. In a forthcoming paper, we shall introduce Hall-Littlewood Laurent symmetric functions which form a basis of $\cd\Lambda_t$, and correspond to the Hall basis of the derived Hall algebra.
	
	As stated in \S\ref{sec:elliptic}, following \cite{BS12}, the derived Hall algebras of elliptic curves have deep connections to DAHAs. Since derived Hall algebras are intrinsic and categorical, it is worth studying the derived Hall algebras of elliptic curves in detail, and making its connection to DAHAs more clearly. 
	
	In a sequel of this paper, we shall construct (extended) derived Hall algebra $\widetilde{\ch}(\cd_m(\ca))$ for any even periodic derived category $\cd_m(\ca)$ of a
	hereditary abelian category $\ca$, where a suitable quantum torus should be appended. Actually, when $m=2$, it coincides with the extended derived Hall algebra as stated  in Remark \ref{rem}. The relationship between  $\widetilde{\ch}(\cd_m(\ca))$ and semi-derived Hall algebra defined in \cite{LP21} shall be studied.

	\subsection{Organization}
	
	The paper is organized as follows.
	In Section \ref{sec:prel}, we recall the definition of root categories, and study some properties related to Hall algebras.  Section \ref{sec:derivedhall} is devoted to defining the derived Hall algebras of root categories, and proving the associativity. In 
	Section \ref{sec:Drinfelddouble}, we prove the derived Hall algebra is isomorphic to the Drinfeld double of Hall algebra. In Sections \ref{sec:Jordan}--\ref{sec:elliptic}, we study the derived Hall algebras of the Jordan quiver and elliptic curves and their connections to symmetric polynomials and DAHAs briefly.

	\subsection{Acknowledgments}
	
	ML deeply thanks Liangang Peng and Weiqiang Wang for guiding him to study representation theory, and also their continuing encouragement.
	ML thanks  Shanghai Key Laboratory of Pure Mathematics and Mathematical Practice, East China Normal University for hospitality
	and support. SR deeply thanks Bangming Deng and Jie Xiao for bringing him to the topic of Hall algebras and DAHAs. SR thanks Haicheng Zhang for careful reading and helpful comments on the original version of the paper. We thank the
	anonymous referee for very helpful suggestions and comments. 
	
	JC is partially supported by the National Natural Science Foundation of China
	(grant No. 11971398). 
	ML is partially supported by the Science and Technology Commission of Shanghai Municipality (grant No. 18dz2271000), and the National Natural Science Foundation of China
	(grant No. 12171333). 
	SR is partially supported by the Natural Science Foundation of Xiamen, China (grant No. 3502Z20227184), 
	the Natural Science Foundation of Fujian Province of China (grant No. 2022J01034), the Fundamental Research Funds for Central Universities of China (grant No. 20720220043) and
	the National Natural Science Foundation of China (grant No. 12271448).

	\section{Preliminaries}
	\label{sec:prel}
	
	For a set $S$, we denote by $|S|$ its cardinality. 
	For any (essentially small) abelian category $\cb$, we denote by $\Iso(\cb)$ the set of isomorphism classes $[M]$ of objects $M\in\cb$. Let $K_0(\cb)$ be the Grothendieck group, and denote by $\widehat{M}$ the class in $K_0(\cb)$ for any object $M\in\cb$. 
	
	\subsection{Root categories}
	Let $\bfk=\mathbb F_q$ be a finite field of $q$ elements. 
	In the following, we always assume that $\ca$ is a hereditary $\bfk$-linear abelian  category which is essentially small with finite-dimensional morphism and extension spaces. 
	
	Let $\langle \cdot,\cdot\rangle$ be the Euler form of $\ca$, i.e., 
	\begin{align}
		\langle M,N\rangle=\frac{|\Hom_\ca(M,N)|}{|\Ext^1_\ca(M,N)|},\quad \forall M,N\in\ca.
	\end{align}
	The Euler form descends to Grothendieck group $K_0(\ca)$, which is also denoted by $\langle\cdot,\cdot\rangle$. We call the Euler form skew symmetric if $\langle M,N\rangle\langle N,M\rangle=1$ for any $M,N\in\ca$.

	Let $\cd^b(\ca)$ be the bounded derived category of $\ca$ with $[1]$ the shift functor. It is well known that $\cd^b(\ca)$
	is a triangulated category. Let $F$ be an automorphism of $\cd^b(\ca)$. The orbit category $\cd^b(\ca)/F$ has the same objects as $\cd^b(\ca)$ and 
	\begin{align*}
		\Hom_{\cd^b(\ca)/F}(X^\bullet,Y^\bullet)=\bigoplus_{i\in\Z}\Hom_{\cd^b(\ca)}(X^\bullet,F^iY^\bullet).
	\end{align*}
	For any $m\ge1$, the orbit category $\cd_m(\ca):=\cd^b(\ca)/[m]$ is called the $m$-periodic derived category of $\ca$, and 
	it is a triangulated category such that the natural projection $\cd^b(\ca)\rightarrow \cd_m(\ca)$ is a triangulated functor; see \cite{PX97,Ke05,St17}. We also denote the shift functor of $\cd_m(\ca)$ by $[1]$. 
	For $m=2$, we also denote $\mathcal{R}(\ca):=\cd_2(\ca)=\cd^b(\ca)/[2]$, which is also called the root category. 
	
	Any object of $\ca$ can be viewed as a stalk complex concentrated at degree zero. This induces a full embedding $\ca\rightarrow \cd^b(\ca)$, and then a full embedding $\ca\rightarrow \cR(\ca)$. We always view $\ca$ as a full subcategory of $\cR(\ca)$ in this way. 
	
	For any object $X^\bullet\in\cR(\ca)$, it is isomorphic to $X_0\oplus X_1[1]$ for some (unique up to isomorphism) $X_0,X_1\in\ca$. In this way, we call $X_i$ the $i$-th cohomology  group of $X^\bullet$, and denote it by $H^i(X^\bullet)$ for $i=0,1$. For a morphism $f:X^\bullet\rightarrow Y^\bullet$ in $\cR(\ca)$, it is of the form 
	\begin{align*}
		f=\begin{bmatrix} f_{00} &f_{01}[1]\\
			f_{10}&f_{11}[1]   
		\end{bmatrix}:X_0\oplus X_1[1]\rightarrow Y_0\oplus Y_1[1]
	\end{align*}
	up to isomorphism, where $f_{00}:X_{0}\rightarrow Y_0$, $f_{01}:X_1\rightarrow Y_0[1]$, $f_{10}:X_0\rightarrow Y_1[1]$ and $f_{11}:X_1\rightarrow Y_1$ are morphisms in $\cd^b(\ca)$. In particular, $f_{00}$ and $f_{11}$ can be viewed as in $\ca$, and we denote by
	$H^0(f)=f_{00}:X_0\rightarrow Y_0$ and $H^1(f)=f_{11}:X_1\rightarrow Y_1$.
	
	We denote by
	$\widehat{X^\bullet}=\widehat{X_0}-\widehat{X_1}\in K_0(\ca)$.
	For any triangle $X^\bullet\rightarrow Y^\bullet\rightarrow Z^\bullet\rightarrow X^\bullet[1]$, we have
	$\widehat{X^\bullet}-\widehat{Y^\bullet}+\widehat{Z^\bullet}=0$. 
	For any $X^\bullet,Y^\bullet\in \cR(\ca)$, we define 
	\begin{align}
		\langle X^\bullet,Y^\bullet\rangle:=\langle \widehat{X^\bullet},\widehat{Y^\bullet}\rangle 
		=&\frac{\langle H^0(X^\bullet),H^0(Y^\bullet)\rangle \langle H^1(X^\bullet),H^1(Y^\bullet)\rangle}{\langle H^0(X^\bullet),H^1(Y^\bullet)\rangle\langle H^1(X^\bullet),H^0(Y^\bullet)\rangle}.
	\end{align}

	\subsection{A useful proposition}
	
	For two objects $X^\bullet$ and $Z^\bullet$ in $\cR(\ca)$, denote 
	\begin{align}
		\{X^\bullet,Z^\bullet\}:=|\Hom_{\ca}(H^0(X^\bullet),H^0(Z^\bullet))|\cdot|\Ext^1_{\ca}(H^1(X^\bullet),H^1(Z^\bullet))|.
	\end{align}
	If $X^\bullet=X_0\oplus X_1[1]$ and $Z^\bullet=Z_0\oplus Z_1[1]$ for $X_0,X_1,Z_0,Z_1\in \ca \subset \mathcal{R}(\ca)$, then 
	\begin{align}
		\{X^\bullet,Z^\bullet\}=|\Hom_{\ca}(X_0,Z_0)|\cdot|\Ext^1_{\ca}(X_1,Z_1)|.
	\end{align}

	We denote by $\Im(X^\bullet,f)$ the image of the map $\Hom_{\cR(\ca)}(X^\bullet,f):\Hom_{\cR(\ca)}(X^\bullet,Y^\bullet)\rightarrow \Hom_{\cR(\ca)}(X^\bullet,Z^\bullet)$ for any morphism $f:Y^\bullet\rightarrow Z^\bullet$ in $\cR(\ca)$. The set $\Im(f,X^\bullet)$ is defined dually.

	\begin{proposition}\label{prop1}
		Given an object $X^\bullet$ and a triangle \begin{align}
			\label{eq:triangle}
			Z^\bullet\stackrel{l}{\longrightarrow}M^\bullet\stackrel{m}{\longrightarrow}L^\bullet\stackrel{n}{\longrightarrow}Z^\bullet[1]
		\end{align}
		in $\mathcal{R}(\ca)$, we have 
		\begin{align}
			\label{eq1}
			|\Im (X^\bullet,n)|=&\dfrac{\{X^\bullet,Z^\bullet\}\{X^\bullet,L^\bullet\}}{\{X^\bullet,M^\bullet\}\langle \widehat{X^\bullet},\delta_0\rangle},
			\\
			|{\mathrm {Im} }(n,X^\bullet)|=&\dfrac{\{Z^\bullet[1],X^\bullet\}\{L^\bullet[1],X^\bullet\}}{\{M^\bullet[1],X^\bullet\}\langle\delta_0, \widehat{X^\bullet}\rangle},
		\end{align}
		where $\delta_0=\widehat{\ker  H^0(l)}$. 
	\end{proposition}	
	
	\begin{proof}
		We only prove  \eqref{eq1}. The second equation can be proved similarly.
		
		For simplicity, denote the functor $\Hom_{\cR(\ca)}(X^\bullet,-)$ by $(X^\bullet,-)$.
		
		We assume $X^\bullet=X_0\oplus X_1[1] $ for $X_0,X_1\in \ca$. Then $\Im (X^\bullet,n)=\Im (X_0,n)\oplus \Im (X_1[1],n)$ and we only need to prove (\ref{eq1}) for $X^\bullet=X_0$ and $X^\bullet=X_1[1]$, respectively. Firstly, we consider the case that $X^\bullet=X_0$.
		
		Step 1. Assume $L^\bullet=L_0 $, $Z^\bullet=Z_0$  $\in \ca\subset \mathcal{R}(\ca)$.
		
		In this case, $M^\bullet$ must be isomorphic to an object $M_0\in\ca\subset \cR(\ca)$, and $0\rightarrow Z_0\xrightarrow{l} M_0\xrightarrow{m} L_0\rightarrow0$ is a short exact sequence in $\ca$. We have $H^0(l)=l$ which is injective, and then $\delta_0=0$. 
		Note that there are no nonzero morphisms from $X_0$ to $Z_0$ factoring through $L_0[1]$ since $\mathcal{A}$ is hereditary. So
		\begin{align*}
			|\Im (X_0,n[1])|=1.
		\end{align*}
		
		Applying $(X^\bullet,-)$ to \eqref{eq:triangle}, we get the long exact sequence\[\begin{tikzpicture}
			\node (0) at (-5.2,0) {$\cdots$};
			\node (5) at (8.4,0) {$\cdots$.};
			\node (-2) at (-3,0) {$(X^\bullet,Z^\bullet)$};
			\node (2) at (0,0) {$(X^\bullet,M^\bullet)$};
			\draw[->] (-2) --node[above ]{\tiny$(X^\bullet,l)$} (2);
			\node (3) at (3,0) {$(X^\bullet,L^\bullet)$};
			\draw[->] (2) --node[above ]{\tiny$(X^\bullet,m)$} (3);
			\node (4) at (2+4,0) {$(X^\bullet,Z^\bullet[1])$};
			\draw[->] (3) --node[above ]{\tiny$(X^\bullet,n)$} (4);
			\draw[->] (0) --node[above ]{} (-2);
			\draw[->] (4) --node[above ]{} (5);
		\end{tikzpicture}\]
		Hence we have \begin{align*}
			|\Im (X_0,l)|=\dfrac{|(X_0,Z_0)|}{|\Im (X_0,n[1])|}= |(X_0,Z_0)|,
		\end{align*} 
		\begin{align} \label{eq:step1}
			|\Im (X_0,m)|=\dfrac{|(X_0,M_0)|}{|\Im (X_0,l)|}= \dfrac{|(X_0,M_0)|}{|(X_0,Z_0)|} ,
		\end{align} 
		and
		\begin{align*}
			|\Im (X_0,n)|=\dfrac{|(X_0,L_0)|}{|\Im (X_0,m)|}= \dfrac{|(X_0,Z_0)||(X_0,L_0)|}{|(X_0,M_0)|}.
		\end{align*}

		Step 2. Assume $l:Z\rightarrow M$ is a morphism in $\mathcal{R}(\ca)$ for $Z,M\in\ca$.
		
		Then $l:Z\rightarrow M$ is a morphism in $\ca$. 
		Since $\ca$ is an abelian category, there exist an object $I\in\ca$, an epimorphism $p:Z\rightarrow I$ and a monomorphism $i:I\rightarrow M$ with $l=i\circ p$. 
		Using the octahedral axiom, we have triangles in $\cR(\ca)$: 
		\begin{align*}
			I\stackrel{i}{\longrightarrow} M\longrightarrow C\longrightarrow I[1],\qquad K\longrightarrow Z\stackrel{p}{\longrightarrow}I\longrightarrow K[1],
		\end{align*}
		and
		\begin{align*}
			Z\stackrel{l}{\longrightarrow} M\longrightarrow K[1]\oplus C\longrightarrow Z[1],
		\end{align*}
		where $K=\ker p$ and $C=\coker i$.
		
		It is known that $\Hom_{\cR(\ca)}(X_0,i)$ is injective. It follows that\[
		|\Im (X_0,l)|=|\Im (X_0,p)|=\dfrac{|(X_0,Z)|}{|(X_0,K)|},
		\]
		where we use the conclusion \eqref{eq:step1} in Step 1.
		Hence 
		\begin{align}
			\label{eq:step2}
			|\Im (X_0,l[1])|=\dfrac{|(X_0,Z)||(X_0,M[1])|}{|(X_0,K\oplus C[1])||\Im (X_0,l)|}= \dfrac{|(X_0,M[1])|}{|(X_0,C[1])|}.
		\end{align}

		Step 3. Assume $l=\begin{pmatrix}
			l_0\\l_1
		\end{pmatrix}:Z\rightarrow M_0\oplus M_1[1]$ is a morphism in $\mathcal{R}(\ca)$ for $Z,M_0,M_1\in\ca$.
		
		In this case, there exists a triangle in $\cR(\ca)$
		\[N\stackrel{p_1}{\longrightarrow} Z\stackrel{l_1}{\longrightarrow} M_1[1] \longrightarrow N[1] \]
		with $p_1$ an epimorphism in $\ca$. Using the octahedral axiom, we have
		\[\xymatrix{N\ar@{=}[r]  \ar[d]_{p_1}& N \ar[d]  \\
			Z\ar[r]^{l_0} \ar[d]_{l_1} & M_0\ar[r] \ar[d] & K[1]\oplus C\ar[r] \ar@{=}[d] &Z[1] \ar[d] 
			\\
			M_1[1] \ar[r] \ar[d]& K'[1]\oplus C' \ar[r] \ar[d] & K[1]\oplus C \ar[r]  & M_1 \\
			N[1] \ar@{=}[r] & N[1]}\]
		%
		%
		%
		%
		where $K=\Ker l_0$, $C=\coker l_0$, $K'=\ker (l_0\circ p_1)$, and $C'=\coker (l_0\circ p_1)\cong \coker l_0=C$ since $p_1$ is surjective. It yields the following triangle, which is isomorphic to \eqref{eq:triangle}, and we identify them in the following: 
		\begin{align*}
			Z\stackrel{l}{\longrightarrow}M_0\oplus M_1[1]\stackrel{m}{\longrightarrow} K'[1]\oplus C'\stackrel{n}{\longrightarrow}Z[1].
		\end{align*}
		
		Now, we can calculate $|\Im (X_0,l[1])|$. Noting that $\Im (X_0,l_1[1])=0$, it follows that\[
		|\Im (X_0,l[1])|=|\Im (X_0,l_0[1])|
		=\dfrac{|(X_0,M_0[1])|}{|(X_0,C[1])|},
		\]where we use the conclusion \eqref{eq:step2} in Step 2.
		Hence we have
		\begin{align}
			\label{eq:step3-1}
			|\Im (X_0,m[1])|=&\dfrac{|(X_0,M_0[1]\oplus M_1)|}{|\Im (X_0,l[1])|}= {|(X_0,M_1)||(X_0,C[1])|}.
		\end{align}
		A simple computation shows that 
		\begin{align*}
			|\Im (X_0,n[1])|=&\dfrac{|(X_0,K'\oplus C'[1])|}{|\Im(X_0,m[1])|}= \dfrac{|(X_0,K')|}{|(X_0,M_1)|},
			\\
			|\Im (X_0,l)|=&\dfrac{|(X_0,Z)|}{|\Im (X_0,n[1])|}= \dfrac{|(X_0,M_1)||(X_0,Z)|} {|(X_0,K')|}.
		\end{align*}
		Then we have
		\begin{align}
			\label{eq:step3-2}
			|\Im (X_0,m)|=&\dfrac{|(X_0,M_0\oplus M_1[1])|}{|\Im (X_0,l)|}
			= \dfrac{|(X_0,M_0\oplus M_1[1])||(X_0,K')|}{|(X_0,M_1)||(X_0,Z)|} .
		\end{align}
		
		Step 4.  Assume $l=\begin{pmatrix}
			l_0\\l_1
		\end{pmatrix}:Z^\bullet\rightarrow M_0\oplus M_1[1]$ is a morphism in $\mathcal{R}(\ca)$ for $Z^\bullet \in {\mathcal{R}}(\ca)$ and $M_0,M_1\in\ca\subset {\mathcal{R}}(\ca)$.
		
		Using the push-out, we have the following commutative diagram
		\begin{equation}
			\label{eq:commdiag-step4}
			\begin{tikzpicture}
				\node (-2) at (-2.5,0) {$Z^\bullet$};
				\node (2) at (0,0) {$M_0$};
				\draw[->] (-2) --node[above ]{\tiny$l_0$} (2);
				\node (3) at (2.5,0) {$N^\bullet$};
				\draw[->] (2) --node[above ]{\tiny$m_0$} (3);
				\node (4) at (3+1+1,0) {$Z^\bullet[1]$};
				\draw[->] (3) --node[above ]{\tiny$n_0$} (4);
				\draw[-] (2.5,-.3) --node[above ]{} (2.5,-1.7);
				\draw[-] (2.6,-.3) --node[above ]{} (2.6,-1.7);
				\node (-21) at (-2.5,-2) {$M_1[1]$};
				\node (21) at (0,-2) {$L^\bullet$};
				\draw[->] (-21) --node[above ]{\tiny$l_2$} (21);
				\node (31) at (2.5,-2) {$N^\bullet$};
				\draw[->] (21) --node[above ]{\tiny$m_2$} (31);
				\node (41) at (3+1+1,-2) {$M_1$};
				\draw[->] (31) --node[above ]{\tiny$n_2$} (41);
				
				\draw[->] (-2) --node[left ]{\tiny$l_1$} (-21);
				\draw[->] (2) --node[left]{\tiny$m_1$} (21);

				\draw[->] (4) --node[right]{\tiny$l_1[1]$} (41);
			\end{tikzpicture}
		\end{equation}
		and a triangle
		\[\begin{tikzpicture}
			\node (-2) at (-3.5,0) {$Z^\bullet$};
			\node (2) at (0-.5,0) {$M_0\oplus M_1[1]$};
			\draw[->] (-2) --node[above ]{\tiny$l$} (2);
			\node (3) at (2.5,0) {$L^\bullet$};
			\draw[->] (2) --node[above ]{\tiny$(-m_1,l_2)$} (3);
			\node (4) at (3+1+1,0) {$Z^\bullet[1]$};
			\draw[->] (3) --node[above ]{\tiny$n_0\circ m_2$} (4);
		\end{tikzpicture}\]
		which is isomorphic to \eqref{eq:triangle}, and we identify them in the following. In particular, $n=n_0\circ m_2$. 
		
		By the conclusion \eqref{eq:step3-2} in Step 3, we get
		\[
		|\Im (X_0,n_0)|
		=\dfrac{|(X_0,H^0(Z^\bullet))||(X_0,H^0(N^\bullet))|}{|(X_0,M_0)|}\cdot \dfrac{|\Ext^1_{\ca}(X_0,H^1(N^\bullet))|}{|(X_0,H^1(N^\bullet))|},
		\]
		and by \eqref{eq:step3-1} in Step 3, we get
		\[
		|\Im (X_0,m_2)|
		=|({X_0,H^0(L^\bullet))|\cdot|(X_0,H^1(N^\bullet)[1]})|.
		\]
		
		By applying $(X^\bullet,-)$ to the commutative diagram \eqref{eq:commdiag-step4}, we observe that the morphism $(X^\bullet,m_2)$ and $(X^\bullet,n_0)$ 
		satisfy that ${\tilde p}\circ {\tilde i}=0$, where  $(\Im (X^\bullet,n_2),{\tilde p})$ is the cokernel of $(X^\bullet,m_2)$ and $(\Im (X^\bullet,m_0),{\tilde i})$ is the kernel of $(X^\bullet,n_0)$. Therefore, by Proposition \ref{cor1}, we obtain
		\begin{align}
			\label{eq:step4-0n}
			|\Im (X_0,n)|&=|\Im (X_0,n_0\circ m_2)|
			\\\notag
			&=\dfrac{|\Im ({X_0},m_2)|\cdot |\Im ({X_0},n_0)|}{|({X_0},N^\bullet)|}
			\\
			\notag&	=\dfrac{|(X_0,H^0(Z^\bullet))||(X_0,H^0(L^\bullet))|}{|(X_0,M_0)|}\cdot \dfrac{|\Ext^1_{\ca}(X_0,H^1(N^\bullet))|}{|(X_0,H^1(N^\bullet))|}
			\\\notag
			&=\dfrac{|(X_0,H^0(Z^\bullet))||(X_0,H^0(L^\bullet))|}{|(X_0,H^0(M^\bullet))|\langle \widehat{X_0},\widehat{H^1(N^\bullet)}\rangle}.
		\end{align}
		It follows that 
		\begin{align}
			\label{eq:step4-1n}
			|\Im (X_1[1],n)|&	=|\Im (X_1,n[1])|
			\\\notag
			&=\dfrac{|(X_1,Z^\bullet)||(X_1,L^\bullet)|}{|(X_1,M^\bullet)||\Im (X_1,n)|}\\
			\notag &=\dfrac{|\Ext^1_{\ca}(X_1,H^1(Z^\bullet))||\Ext^1_{\ca}(X_1,H^1(L^\bullet))|}{|\Ext^1_{\ca}(X_1,H^1(M^\bullet))|}\cdot \langle \widehat{X_1},\widehat{H^1(N^\bullet)}\rangle.
		\end{align}
		Combining \eqref{eq:step4-0n} with \eqref{eq:step4-1n}, we have
		\[
		|\Im (X^\bullet,n)|=\dfrac{\{X^\bullet,Z^\bullet\}\{X^\bullet,L^\bullet\}}{\{X^\bullet,M^\bullet\}\langle \widehat{X^\bullet},\widehat{H^1(N^\bullet)}\rangle}.
		\]
		
		It remains to prove that $\widehat{H^1(N^\bullet)}=\delta_0$. 
		In fact, from \eqref{eq:triangle} and the first row in \eqref{eq:commdiag-step4}, we have the long exact sequences
		\[\begin{tikzpicture}
			\node (-3) at (-4*1.5+1,0) {$\cdots$};
			\node (-2) at (-2*1.5,0) {$H^1(L^\bullet)$};
			\node (2) at (0,0) {$H^0(Z^\bullet)$};
			\draw[->] (-2) --node[above ]{} (2);
			\node (3) at (2*1.5,0) {$H^0(M^\bullet)$};
			\draw[->] (2) --node[above ]{\tiny$H^0(l)$} (3);
			\node (4) at (6,0) {$H^0(L^\bullet)$};
			\draw[->] (3) --node[above ]{} (4);
			\node (5) at (8,0) {$\cdots$,};
			\draw[->] (4) --node[above ]{} (5);
			\draw[->] (-3) --node[above ]{} (-2);
		\end{tikzpicture}\]and
		\[\begin{tikzpicture}
			\node (-3) at (-5,0) {$\cdots$};
			\node (-2) at (-3,0) {$0$};
			\node (2) at (0,0) {$H^1(N^\bullet)$};
			\draw[->] (-2) --node[above ]{} (2);
			\node (3) at (3,0) {$H^0(Z^\bullet)$};
			\draw[->] (2) --node[above ]{} (3);
			\node (4) at (6,0) {$H^0(M_0)$};
			\draw[->] (3) --node[above ]{\tiny$H^0(l_0)$} (4);
			\node (5) at (8,0) {$\cdots$.};
			\draw[->] (4) --node[above ]{} (5);
			\draw[->] (-3) --node[above ]{} (-2);
		\end{tikzpicture}\]
		Hence $\delta_0=\widehat{\ker  H^0(l)}=\widehat{\ker  H^0(l_0)}=\widehat{H^1(N^\bullet)}$.
		
		The proof is completed.
	\end{proof}

	\section{Derived Hall algebras of root categories}
	\label{sec:derivedhall}
	
	In this section, unless otherwise specific, we always assume that the Euler form of $\ca$ is \emph{skew symmetric}. We shall define a kind of derived Hall algebras for the root category of $\ca$ in this section.

	\subsection{Derived Hall numbers}
	
	We recall some notations in the root category $\cR(\ca)$ from \cite{PX00,XX08}. 
	
	Let $Z^\bullet,M^\bullet\in\cR(\ca)$. Define the radical of $\Hom(Z^\bullet,M^\bullet)$ to be 
	\begin{align*}
		\rad\Hom(Z^\bullet,M^\bullet):=\{
		l\in \Hom(Z^\bullet,M^\bullet)\mid g\circ l\circ h \text{ is not an isomorphism for any }
		\\
		h:X^\bullet\rightarrow Z^\bullet \text{ and }g:M^\bullet\rightarrow X^\bullet\text{ with }X^\bullet\in \cR(\ca) \text{ indecomposable}\}.  
	\end{align*}
	For any $M^\bullet,L^\bullet,Z^\bullet\in\cR(\ca)$, and $\delta,\delta'\in K_0(\ca)$, we define
	\begin{align}
		\label{eq:HomZMLdelta}
		\cw(Z^\bullet,L^\bullet;M^\bullet)_\delta&=\{(l,m,n)\in\Hom(Z^\bullet,M^\bullet)\times \Hom(M^\bullet,L^\bullet)\times\Hom(L^\bullet,Z^\bullet[1])\mid 
		\\\notag
		&Z^\bullet \xrightarrow{l}M^\bullet\xrightarrow{m} L^\bullet\xrightarrow{n}Z^\bullet[1] \text{ is a triangle and } \widehat{\coker  H^0(m)}=\delta \}.
	\end{align}
	Let $\Aut(Z^\bullet)$ be the automorphism group of $Z^\bullet$. Then there is a natural action of $\Aut(Z^\bullet)\times \Aut(L^\bullet)$ on $\cw(Z^\bullet,L^\bullet;M^\bullet)_\delta$, 
	and its orbit space is denoted by
	\begin{align}
		\cv(Z^\bullet,L^\bullet;M^\bullet)_\delta=&\cw(Z^\bullet,L^\bullet;M^\bullet)_\delta/(\Aut (Z^\bullet)\times \Aut(L^\bullet)).
	\end{align}
	For convenience, we also denote 
	\begin{align}
		\cw'(Z^\bullet,L^\bullet;M^\bullet)_{\delta'}&=\{(l,m,n)\in\Hom(Z^\bullet,M^\bullet)\times \Hom(M^\bullet,L^\bullet)\times\Hom(L^\bullet,Z^\bullet[1])\mid 
		\\\notag
		&Z^\bullet \xrightarrow{l}M^\bullet\xrightarrow{m} L^\bullet\xrightarrow{n}Z^\bullet[1] \text{ is a triangle and } \widehat{\coker  H^0(l)}=\delta' \},
	\end{align}
	which is the same as $\cw(L^\bullet[1],M^\bullet;Z^\bullet)_{\delta'}$, 
	and 
	$\cv'(Z^\bullet,L^\bullet;M^\bullet)_{\delta'}$ 
	its orbit space under the action of $\Aut (Z^\bullet)\times \Aut(L^\bullet)$. Also define
	\begin{align}
		\Hom(Z^\bullet,M^\bullet)_{L^\bullet,\delta'}&=\{l:Z^\bullet\rightarrow M^\bullet\mid \cone l\cong L^\bullet, \widehat{\coker  H^0(l)}=\delta'\}.
	\end{align}

	By \cite[Proposition 2.6]{XX08}, whose proof holds for arbitrary triangulated categories, we obtain the following bijections
	\begin{align}
		\label{eq:bijection1}
		\Hom(M^\bullet,L^\bullet)_{Z^\bullet[1],\delta}&\longleftrightarrow
		\cw(Z^\bullet,L^\bullet;M^\bullet)_\delta/ \Aut(Z^\bullet),
		\\
		\Hom(M^\bullet,L^\bullet)_{Z^\bullet[1],\delta} / \Aut(L^\bullet)
		&\longleftrightarrow\cv(Z^\bullet,L^\bullet;M^\bullet)_\delta,
		\\
		\Hom(Z^\bullet,M^\bullet)_{L^\bullet,\delta'}&\longleftrightarrow
		\cw'(Z^\bullet,L^\bullet;M^\bullet)_{\delta'}	/ \Aut(L^\bullet),
		\\
		\Hom(Z^\bullet,M^\bullet)_{L^\bullet,\delta'}	/ \Aut(Z^\bullet)
		&\longleftrightarrow	\cv'(Z^\bullet,L^\bullet;M^\bullet)_{\delta'}.
	\end{align}
	
	\begin{lemma}
		\label{lem:W=W'}
		For any $M^\bullet,L^\bullet,Z^\bullet\in\cR(\ca)$, and $\delta,\delta'\in K_0(\ca)$, if $\delta+\delta'=\widehat{H^0(L^\bullet)}$, then 
		\begin{align*}
			\cw(Z^\bullet,L^\bullet;M^\bullet)_{\delta}=	\cw'(Z^\bullet,L^\bullet;M^\bullet)_{\delta'},
		\end{align*}
		and hence
		\begin{align*}
			\cv(Z^\bullet,L^\bullet;M^\bullet)_{\delta}=	\cv'(Z^\bullet,L^\bullet;M^\bullet)_{\delta'}.
		\end{align*}
	\end{lemma}
	
	\begin{proof}
		First we prove that $\cw(Z^\bullet,L^\bullet;M^\bullet)_\delta\subseteq \cw'(Z^\bullet,L^\bullet;M^\bullet)_{\delta'}$. In fact, for any $(l,m,n)\in \cw(Z^\bullet,L^\bullet;M^\bullet)_\delta$, we have $\widehat{\coker  H^0(l)}=\widehat{L_0}-\widehat{\coker  H^0(m)}=\widehat{L_0}-\delta=\delta'$. Therefore $\cw(Z^\bullet,L^\bullet;M^\bullet)_\delta\subseteq \cw'(Z^\bullet,L^\bullet;M^\bullet)_{\delta'}$. By a similar argument, we have $ \cw'(Z^\bullet,L^\bullet;M^\bullet)_{\delta'}\subseteq
		\cw(Z^\bullet,L^\bullet;M^\bullet)_\delta$ and hence the first statement follows. 
		
		By the action of $\Aut (Z^\bullet)\times \Aut(L^\bullet)$, we get the second statement.
	\end{proof}
	
	In order to introduce the derived Hall numbers,  we recall from \cite{XX08} the following decomposition of a triangle, whose proof also holds for root categories. 
	
	\begin{lemma}
		[\text{\cite[Remark 2.3]{XX08}}]
		\label{lem:representative}
		For any $M^\bullet,L^\bullet,Z^\bullet\in \cR(\ca)$ and any $\alpha\in \cv(Z^\bullet,L^\bullet;M^\bullet)$, we have a representative of $\alpha$ with the following form:
		\[
		\begin{tikzpicture}
			\node (-2) at (-3,0) {$Z^\bullet$};
			\node (2) at (0,0) {$M^\bullet$};
			\draw[->] (-2) --node[above ]{$(0,l_2)$} (2);
			\node (3) at (3,0) {$L^\bullet$};
			\draw[->] (2) --node[above ]{\tiny$\begin{pmatrix}
					0\\m_2
				\end{pmatrix}$} (3);
			\node (4) at (6,0) {$Z^\bullet[1]$,};
			\draw[->] (3) --node[above ]{\tiny$\begin{pmatrix}
					n_{11}&0\\0&n_{22}
				\end{pmatrix}$} (4);
		\end{tikzpicture}
		\] where $Z^\bullet\cong Z_1^\bullet(\alpha)\oplus Z_2^\bullet(\alpha)$, $L^\bullet\cong L_1^\bullet(\alpha)\oplus L_2^\bullet(\alpha)$, $n_{11}:L^\bullet_1(\alpha)\rightarrow Z^\bullet_1(\alpha)[1]$ is an isomorphism  and $n_{22}\in \mathrm{rad} \Hom(L^\bullet_2(\alpha),Z^\bullet_2(\alpha)[1])$.
	\end{lemma}

	With the help of Lemma \ref{lem:representative}, we can state the following proposition, which is an analog of \cite[Proposition 2.5']{XX08}.
	
	\begin{proposition}
		\label{prop4.1}
		For any $M^\bullet,L^\bullet,Z^\bullet\in\cR(\ca)$ and $\delta,\delta'\in K_0(\ca)$, we have
		\begin{align*}
			\dfrac{|\Hom(M^\bullet,L^\bullet)_{Z^\bullet[1],\delta}|}{\langle \widehat{H^0(Z^\bullet)}+\widehat{H^0(L^\bullet)}-\widehat{H^0(M^\bullet)}-\delta,\widehat{L^\bullet}\rangle}\dfrac{\{L^\bullet[1],L^\bullet\}  \{Z^\bullet[1],L^\bullet\} }{|\Aut(L^\bullet)|\{M^\bullet[1],L^\bullet\}}=\sum_{\alpha\in	\cv(Z^\bullet,L^\bullet;M^\bullet)_\delta}&\dfrac{|\End(L^\bullet_1(\alpha))|}{|\Aut(L^\bullet_1(\alpha))|},
			\\
			\dfrac{|\Hom(Z^\bullet,M^\bullet)_{L^\bullet,\delta'}|}{\langle \widehat{Z^\bullet[1]}, \widehat{H^0(Z^\bullet)}-\widehat{H^0(M^\bullet)}+\delta'\rangle}\dfrac{\{Z^\bullet[1],Z^\bullet\}  \{Z^\bullet[1],L^\bullet\} }{|\Aut(Z^\bullet)|\{Z^\bullet[1],M^\bullet\}}=\sum_{\alpha\in	\cv'(Z^\bullet,L^\bullet;M^\bullet)_{\delta'}}&\dfrac{|\End(L^\bullet_1(\alpha))|}{|\Aut(L^\bullet_1(\alpha))|},
		\end{align*}
		where $(l,m,n)$ is a representative of $\alpha$ such that $n:L^\bullet\rightarrow Z^\bullet[1]$ is of the form
		$$
		n=\begin{pmatrix}
			n_{11}&0\\0&n_{22}
		\end{pmatrix}:L^\bullet_1(\alpha)\oplus L^\bullet_2(\alpha)\longrightarrow Z^\bullet_1(\alpha)[1]\oplus Z^\bullet_2(\alpha)[1]
		$$ 
		with
		$n_{11}:L^\bullet_1(\alpha)\rightarrow Z^\bullet_1(\alpha)[1]$ an isomorphism  and $n_{22}\in \mathrm{rad} \Hom(L^\bullet_2(\alpha),Z^\bullet_2(\alpha)[1])$.
	\end{proposition}
	
	\begin{proof}
		From \eqref{eq:bijection1},  there exists a bijection
		\[
		\Hom(M^\bullet,L^\bullet)_{Z^\bullet[1],\delta}\longrightarrow
		\cw(Z^\bullet,L^\bullet;M^\bullet)_\delta/ \Aut(Z^\bullet).
		\] 
		For any $(l,m,n)\in\cw(Z^\bullet,L^\bullet;M^\bullet)_\delta$, we denote by $(l,m,n)^*$ its orbit under the action of $\Aut(Z^\bullet)$.
		Consider the action
		\begin{align*}
			\Aut(L^\bullet)\times \big(\cw(Z^\bullet,L^\bullet;M^\bullet)_\delta/\Aut(Z^\bullet)\big)&\longrightarrow \cw(Z^\bullet,L^\bullet;M^\bullet)_\delta/\Aut(Z^\bullet)
			\\
			(b,(l,m,n)^*)&\mapsto (l,b^{-1}m,nb)^*,
		\end{align*}
		whose orbit space is just $V(Z^\bullet,L^\bullet;M^\bullet)_\delta$. 
		
		For $(l,m,n)^*\in \cw(Z^\bullet,L^\bullet;M^\bullet)_\delta/\Aut(Z^\bullet)$, similar to the argument in \cite[Proposition 2.5]{XX08}, its stable subgroup, denoted by $G((l,m,n)^*)$, satisfies
		\[
		|G((l,m,n)^*)|=|\Im(n,L^\bullet)|\cdot\dfrac{|\Aut(L^\bullet_1(\alpha))|}{|\End(L^\bullet_1(\alpha))|}.
		\]
		Hence 
		\begin{align*}
			{|\Hom(M^\bullet,L^\bullet)_{Z^\bullet[1],\delta}|}
			&=
			|{W(Z^\bullet,L^\bullet;M^\bullet)_\delta}/{\Aut(Z^\bullet)}|
			\\
			&=\sum_{\alpha\in	\cv(Z^\bullet,L^\bullet;M^\bullet)_\delta}\dfrac{|\Aut(L^\bullet)||\End(L^\bullet_1(\alpha))|}{|\Im(n,L^\bullet)||\Aut(L^\bullet_1(\alpha))|}.
		\end{align*}
		Note that $\delta=\widehat{\coker  H^0(m)}=\widehat{H^0(Z^\bullet)}+\widehat{H^0(L^\bullet)}-\widehat{H^0(M^\bullet)}-\widehat{\ker  H^0(l)}$. Then the first statement follows by Proposition \ref{prop1}. 
		
		The second statement is similar.
	\end{proof}
	
	For any $L^\bullet\cong L_0\oplus L_1[1]$, $M^\bullet\cong M_0\oplus M_1[1]$, $Z^\bullet\cong Z_0\oplus Z_1[1]$ in $\cR(\ca)$ with $L_i,M_i,Z_i\in\ca$ for $i=0,1$, similar to \cite{T06,XX08},
	we define 
	\begin{align}
		\label{3.9}
		F_{L^\bullet Z^\bullet}^{M^\bullet}:&=\dfrac{1}{\{Z^\bullet,L^\bullet[1]\}}
		\sum_{\delta}\ \langle-\delta,\widehat{M^\bullet} \rangle \sum_{\alpha\in	\cv(Z^\bullet,L^\bullet;M^\bullet)_\delta}\dfrac{|\End(L^\bullet_1(\alpha))|}{|\Aut(L^\bullet_1(\alpha))|}.
	\end{align}
	
	Note that 
	$$\dfrac{\{Z^\bullet,L^\bullet[1]\}}{\{Z^\bullet[1],L^\bullet\}}=\dfrac{|\Hom_{\ca}(Z_0,L_1)|\cdot|\Ext^1_{\ca}(Z_1,L_0)|}{|\Hom_{\ca}(Z_1,L_0)|\cdot|\Ext^1_{\ca}(Z_0,L_1)|}=\dfrac{\langle \widehat{Z_0},\widehat{L_1} \rangle}{\langle \widehat{Z_1},\widehat{L_0} \rangle},$$
	and the Euler form of $\ca$ is {skew symmetric}. 
	It follows that 
	\begin{align*}\dfrac{1}{\{Z^\bullet,L^\bullet[1]\}}
		=
		\dfrac{\langle \widehat{Z_0},\widehat{L_0}\rangle\langle \widehat{L_0},\widehat{Z_0}\rangle }{\{Z^\bullet,L^\bullet[1]\}}
		=
		\dfrac{\langle \widehat{Z_0},\widehat{L_1} \rangle}{\langle \widehat{Z_1},\widehat{L_0} \rangle}\cdot\dfrac{\langle \widehat{Z_0},\widehat{L^\bullet}\rangle\langle \widehat{L_0},\widehat{Z^\bullet}\rangle }{\{Z^\bullet,L^\bullet[1]\}}
		=\dfrac{\langle \widehat{Z_0},\widehat{L^\bullet}\rangle\langle \widehat{L_0},\widehat{Z^\bullet}\rangle }{\{Z^\bullet[1],L^\bullet\}}.
	\end{align*}
	Hence formula (\ref{3.9}) can be expressed as
	\begin{align*}
		F_{L^\bullet Z^\bullet}^{M^\bullet}=
		\dfrac{\langle \widehat{Z_0},\widehat{L^\bullet}\rangle\langle \widehat{L_0},\widehat{Z^\bullet}\rangle }{\{Z^\bullet[1],L^\bullet\}}\sum_{\delta}\quad\langle-\delta,\widehat{M^\bullet} \rangle\sum_{\alpha\in	\cv(Z^\bullet,L^\bullet;M^\bullet)_\delta}\dfrac{|\End(L^\bullet_1(\alpha))|}{|\Aut(L^\bullet_1(\alpha))|}.
	\end{align*}

	By Proposition \ref{prop4.1} and Lemma \ref{lem:W=W'}, we have
	\begin{align}
		\label{3.10}
		F_{L^\bullet,Z^\bullet}^{M^\bullet}
		&={\langle \widehat{Z_0},\widehat{L^\bullet}\rangle\langle \widehat{L_0},\widehat{Z^\bullet}\rangle }\sum_{\delta}\langle-\delta,\widehat{M^\bullet} \rangle
		\dfrac{|\Hom(M^\bullet,L^\bullet)_{Z^\bullet[1],\delta}|}{\langle \widehat{Z_0}+\widehat{L_0}-\widehat{M_0}-\delta,\widehat{L^\bullet}\rangle}\cdot \dfrac{\{L^\bullet[1],L^\bullet\}   }{|\Aut(L^\bullet)|\{M^\bullet[1],L^\bullet\}}
		\\\label{3.11}
		&={\langle \widehat{Z_0},\widehat{L^\bullet}\rangle\langle \widehat{L_0},\widehat{Z^\bullet}\rangle }\sum_{\delta'}\langle\delta'-\widehat{L_0},\widehat{M^\bullet} \rangle\dfrac{|\Hom(Z^\bullet,M^\bullet)_{L^\bullet,\delta'}|}{\langle \widehat{Z^\bullet[1]}, \widehat{Z_0}-\widehat{M_0}+\delta'\rangle}\cdot \dfrac{\{Z^\bullet[1],Z^\bullet\}   }{|\Aut(Z^\bullet)|\{Z^\bullet[1],M^\bullet\}}.	
	\end{align}

	\subsection{Derived Hall algebras}
	
	For any $X^\bullet\in\cR(\ca)$, denote by $[X^\bullet]$ its isomorphism class. Let $\och(\cR(\ca))$ be the $\Q$-space  with the basis $\{u_{[X^\bullet]}\mid X^\bullet\in\mathcal{R}(\ca)\}$.  We define 
	\begin{align}
		u_{[Y^\bullet]}\cdot u_{[X^\bullet]}=\sum_{[L^\bullet]}F_{Y^\bullet X^\bullet}^{L^\bullet}\cdot u_{[L^\bullet]}.
	\end{align}
	
	In this section, we shall prove that $\och(\cR(\ca))$ is an associative algebra with $u_{[0]}$ the unit under the assumption that the Euler form of $\ca$ is skew symmetric. Then $\och(\cR(\ca))$ is called the derived Hall algebra of $\ca$.  
	
	To prove the associativity, we need more notations.
	
	Let $L^\bullet,L'^\bullet, M^\bullet,X^\bullet,Y^\bullet,Z^\bullet\in\cR(\ca)$. 
	
	For  $\delta,\delta_1,\delta_2\in K_0(\ca)$, we define $\Hom(M^\bullet\oplus X^\bullet,L^\bullet)^{Y^\bullet,\delta_1,Z^\bullet[1],\delta_2}_{L'^\bullet[1],\delta}$ to be 
	\begin{align*}
		\big\{(m,f):M^\bullet\oplus X^\bullet\rightarrow L^\bullet\mid &\cone (m,f)\cong L'^\bullet[1],\widehat{\coker  H^0((m,f))}=\delta,
		\cone f\cong Y^\bullet,
		\\
		&  \widehat{\coker  H^0(f)}=\delta_1,
		\cone m\cong Z^\bullet[1], \widehat{\coker  H^0(m)}=\delta_2
		\big\},
	\end{align*}
	and  $\cw(L'^\bullet,L^\bullet;M^\bullet\oplus X^\bullet)^{Y^\bullet,\delta_1,Z^\bullet[1],\delta_2}_\delta$ to be
	\begin{align*}
		\big\{\big(\begin{pmatrix}
			f'\\ -m'
		\end{pmatrix},&(m,f),\theta\big)\mid L'^\bullet\stackrel{\tiny \begin{pmatrix}
				f'\\ -m'
		\end{pmatrix}}{\longrightarrow} M^\bullet\oplus X^\bullet \stackrel{(m,f)}{\longrightarrow} L^\bullet\stackrel{\theta}{\longrightarrow} L'^\bullet[1] \text{ is a triangle, } \cone f\cong Y^\bullet,
		\\
		&\widehat{\coker  H^0((m,f))}=\delta,\widehat{\coker  H^0(f)}=\delta_1,
		\cone m\cong Z^\bullet[1], \widehat{\coker  H^0(m)}=\delta_2 \big\}.
	\end{align*}
	Then there is a natural action of $\Aut (L'^\bullet)\times \Aut(L^\bullet)$ on $\cw(L'^\bullet,L^\bullet;M^\bullet\oplus X^\bullet)^{Y^\bullet,\delta_1,Z^\bullet[1],\delta_2}_\delta$, and we denote its orbit space by $\cv(L'^\bullet,L^\bullet;M^\bullet\oplus X^\bullet)^{Y^\bullet,\delta_1,Z^\bullet[1],\delta_2}_\delta$.
	
	
	Furthermore, for  $\delta',\delta'_1,\delta'_2\in K_0(\ca)$, we also define $\Hom(L'^\bullet,M^\bullet\oplus X^\bullet)^{Y^\bullet,\delta'_1,Z^\bullet[1],\delta'_2}_{L^\bullet,\delta'}$ to be
	\begin{align*}
		\big\{\begin{pmatrix}
			f'\\ -m'
		\end{pmatrix}:L'^\bullet\rightarrow M^\bullet\oplus X^\bullet\mid &\cone \begin{pmatrix}
			f'\\ -m'
		\end{pmatrix}\cong L^\bullet, \widehat{\coker  H^0}(\begin{pmatrix}
			f'\\ -m'
		\end{pmatrix})=\delta',\cone f'\cong Y^\bullet, 
		\\
		&
		\widehat{\coker  H^0(f')}=\delta'_1,
		\cone m'\cong Z^\bullet[1], \widehat{\coker  H^0(m')}=\delta'_2
		\big\},
	\end{align*}
	and set $\cw'(L'^\bullet,L^\bullet;M^\bullet\oplus X^\bullet)^{Y^\bullet,\delta'_1,Z^\bullet[1],\delta'_2}_{\delta'}$ to be 
	\begin{align*}
		\{\big(\begin{pmatrix}
			f'\\ -m'
		\end{pmatrix}&,(m,f),\theta\big)\mid L'^\bullet\stackrel{ \tiny\begin{pmatrix}
				f'\\ -m'
		\end{pmatrix}}{\longrightarrow} M^\bullet\oplus X^\bullet \stackrel{(m,f)}{\longrightarrow} L^\bullet\stackrel{\theta}{\longrightarrow} L'^\bullet[1] \text{ is a triangle, } \cone f'\cong Y^\bullet,
		\\
		& \widehat{\coker  H^0}(\begin{pmatrix}
			f'\\ -m'
		\end{pmatrix})=\delta',\widehat{\coker  H^0(f')}=\delta'_1,
		\cone m'\cong Z^\bullet[1], \widehat{\coker  H^0(m')}=\delta'_2 \}.
	\end{align*}
	For the natural action of $\Aut (L'^\bullet)\times \Aut(L^\bullet)$ on $\cw'(L'^\bullet,L^\bullet;M^\bullet\oplus X^\bullet)^{Y^\bullet,\delta'_1,Z^\bullet[1],\delta'_2}_{\delta'}$, its orbit space is denoted by $\cv'(L'^\bullet,L^\bullet;M^\bullet\oplus X^\bullet)^{Y^\bullet,\delta'_1,Z^\bullet[1],\delta'_2}_{\delta'}$.

	\begin{proposition}
		\label{prop2}
		Keep the notations as above. 
		If $\delta'=\widehat{H^0(L^\bullet)}-\delta, \delta'_1=\delta_1-\delta$ and $\delta'_2=\delta_2-\delta$, then
		\[
		\cw(L'^\bullet,L^\bullet;M^\bullet\oplus X^\bullet)^{Y^\bullet,\delta_1,Z^\bullet[1],\delta_2}_\delta=\cw'(L'^\bullet,L^\bullet;M^\bullet\oplus X^\bullet)^{Y^\bullet,\delta'_1,Z^\bullet[1],\delta'_2}_{\delta'},
		\]
		and hence
		\[
		\cv(L'^\bullet,L^\bullet;M^\bullet\oplus X^\bullet)^{Y^\bullet,\delta_1,Z^\bullet[1],\delta_2}_\delta=\cv'(L'^\bullet,L^\bullet;M^\bullet\oplus X^\bullet)^{Y^\bullet,\delta'_1,Z^\bullet[1],\delta'_2}_{\delta'}.
		\]
	\end{proposition}
	
	\begin{proof}
		The proof is inspired by \cite[Proposition 3.2]{XX08}. 
		
		For any $\big(\begin{pmatrix}
			f'\\ -m'
		\end{pmatrix},(m,f),\theta\big)\in \cw(L'^\bullet,L^\bullet;M^\bullet\oplus X^\bullet)^{Y^\bullet,\delta_1,Z^\bullet[1],\delta_2}_\delta$,
		it can be expressed by the following diagram: 
		\[\xymatrix{& Z^\bullet \ar[d]^l\\
			L'^{\bullet} \ar[r]^{f'} \ar[d]_{m'}& M^\bullet \ar[d]^m \\
			X^\bullet \ar[r] ^f & L^\bullet \ar[r]^g\ar[d]^n & Y^\bullet \ar[r]^h & X^\bullet[1] 
			\\
			&Z^\bullet[1]}\]

		
		
		
		
		
		By using the octahedral axiom and the push-out property of triangulated categories, the above diagram is completed to following commutative diagram: 
		\begin{equation}
			\label{eq:octahedral}
			\xymatrix{Z^\bullet\ar[d]_{l'} \ar@{=}[r]& Z^\bullet \ar[d]^l\\
				L'^{\bullet} \ar[r]^{f'} \ar[d]_{m'}& M^\bullet \ar[d]^m \ar[r]^{g'}& Y^\bullet\ar[r]^{h'}\ar@{=}[d]& L'^\bullet[1]\ar[d]^{m'[1]}\\
				X^\bullet \ar[r] ^f \ar[d]_{n'}& L^\bullet \ar[r]^g\ar[d]^n & Y^\bullet \ar[r]^h & X^\bullet[1] 
				\\
				Z^\bullet[1]\ar@{=}[r]&Z^\bullet[1]}
		\end{equation}
		Note that $\theta=-h'\circ g$. 

		It is clear that $\cone f'\cong\cone f\cong Y^\bullet$ and $\cone m'\cong\cone m\cong Z^\bullet[1]$. Obviously, we obtain that $\widehat{\coker  H^0}(\begin{pmatrix} 	f'\\ -m' 	\end{pmatrix})=\widehat{L_0}-\widehat{\coker  H^0((m,f))}=\widehat{L_0}-\delta=\delta'$. Using Corollary \ref{cor1}, we have
		\begin{align*}
			(\widehat{Y_0}-\widehat{\coker  H^0(f')})+\widehat{\coker  H^0(f)}&=\widehat{Y_0}+\widehat{\coker  H^0((m,f))},
			\\
			(\widehat{Z_1}-\widehat{\coker  H^0(m')})+\widehat{\coker  H^0(m)}&=\widehat{Z_1}+\widehat{\coker  H^0((m,f))}.
		\end{align*}
		It follows that 
		\begin{align*}
			\widehat{\coker  H^0(f')}&=\widehat{\coker  H^0(f)}-\widehat{\coker  H^0((m,f))}={\delta_1-\delta}=\delta'_1,
			\\
			\widehat{\coker  H^0(m')}&=\widehat{\coker  H^0(m)}-\widehat{\coker  H^0((m,f))}={\delta_2-\delta}=\delta'_2.
		\end{align*}
		Therefore, 
		\begin{align*}
			\cw(L'^\bullet,L^\bullet;M^\bullet\oplus X^\bullet)^{Y^\bullet,\delta_1,Z^\bullet[1],\delta_2}_\delta\subseteq \cw'(L'^\bullet,L^\bullet;M^\bullet\oplus X^\bullet)^{Y^\bullet,\delta'_1,Z^\bullet[1],\delta'_2}_{\delta'}.
		\end{align*}
		By a similar argument, one can prove have 
		\begin{align*} 
			\cw'(L'^\bullet,L^\bullet;M^\bullet\oplus X^\bullet)^{Y^\bullet,\delta'_1,Z^\bullet[1],\delta'_2}_{\delta'}\subseteq
			\cw(L'^\bullet,L^\bullet;M^\bullet\oplus X^\bullet)^{Y^\bullet,\delta_1,Z^\bullet[1],\delta_2}_\delta
		\end{align*}
		and then the first statement follows.

		By the action of $\Aut (L'^\bullet)\times \Aut(L^\bullet)$, we get the second statement.
	\end{proof}

	\begin{proposition}
		\label{prop:Hallnumbercounting}
		Keep the notations as above. 
		We have
		\begin{align*}
			\dfrac{|\Hom(M^\bullet\oplus X^\bullet,L^\bullet)^{Y^\bullet,\delta_1,Z^\bullet[1],\delta_2}_{L'^\bullet[1],\delta}|}{\langle \widehat{L'_0}+\widehat{L_0}-\widehat{M_0}-\widehat{X_0}-\delta,\widehat{L^\bullet}\rangle}\cdot \dfrac{\{L^\bullet[1],L^\bullet\}  \{L'^\bullet[1],L^\bullet\} }{|\Aut(L^\bullet)|\{M^\bullet[1]\oplus X^\bullet[1],L^\bullet\}}&=\sum_{\alpha\in \cv	}\dfrac{|\End(L^\bullet_1(\alpha))|}{|\Aut(L^\bullet_1(\alpha))|},
			\\
			\dfrac{|\Hom(L'^\bullet,M^\bullet\oplus X^\bullet)^{Y^\bullet,\delta'_1,Z^\bullet[1],\delta'_2}_{L^\bullet,\delta'}|}{\langle \widehat{L'^\bullet[1]}, \widehat{L'_0}-\widehat{M_0}-\widehat{X_0}+\delta'\rangle}\cdot \dfrac{\{L'^\bullet[1],L'^\bullet\}  \{L'^\bullet[1],L^\bullet\} }{|\Aut(L'^\bullet)|\{L'^\bullet[1],M^\bullet\oplus X^\bullet\}}&=\sum_{\alpha\in \cv'}\dfrac{|\End(L^\bullet_1(\alpha))|}{|\Aut(L^\bullet_1(\alpha))|},
		\end{align*}
		where $\cv=\cv(L'^\bullet,L^\bullet;M^\bullet\oplus X^\bullet)^{Y^\bullet,\delta_1,Z^\bullet[1],\delta_2}_\delta$, $\cv'=\cv'(L'^\bullet,L^\bullet;M^\bullet\oplus X^\bullet)^{Y^\bullet,\delta'_1,Z^\bullet[1],\delta'_2}_{\delta'}$, and the triangle $\big(\begin{pmatrix} 	f'\\ -m' 	\end{pmatrix},(m,f),\theta\big)$ is a representative of $\alpha$ such that $\theta:L^\bullet\rightarrow L'^{\bullet}$ is of the form  
		$$
		\theta=\begin{pmatrix}
			\theta_{11}&0\\0&\theta_{22}
		\end{pmatrix}:L^\bullet_1(\alpha)\oplus L^\bullet_2(\alpha)\longrightarrow L'^\bullet_1[1](\alpha)\oplus L'^\bullet_2[1](\alpha)
		$$ 
		with
		$\theta_{11}:L^\bullet_1(\alpha) \to  L'^\bullet_1(\alpha)[1]$ an isomorphism and $\theta_{22}\in \mathrm{rad} \Hom(L^\bullet_2(\alpha),L'^\bullet_2(\alpha)[1])$.
	\end{proposition}
	
	\begin{proof}
		The proof is the same as that of Proposition \ref{prop4.1}, we omit it here.
	\end{proof}
	
	Using Propositions \ref{prop2}--\ref{prop:Hallnumbercounting}, we have the following corollary immediately.
	\begin{corollary}
		\label{cor2}
		If $\delta'=\widehat{L_0}-\delta, \delta'_1=\delta_1-\delta$ and $\delta'_2=\delta_2-\delta$, then
		\begin{align*}
			&\hspace{-2cm}\dfrac{|\Hom(M^\bullet\oplus X^\bullet,L^\bullet)^{Y^\bullet,\delta_1,Z^\bullet[1],\delta_2}_{L'^\bullet[1],\delta}|}{\langle \widehat{L'_0}+\widehat{L_0}-\widehat{M_0}-\widehat{X_0}-\delta,\widehat{L^\bullet}\rangle}\cdot \dfrac{\{L^\bullet[1],L^\bullet\}  \{L'^\bullet[1],L^\bullet\} }{|\Aut(L^\bullet)|\{M^\bullet[1]\oplus X^\bullet[1],L^\bullet\}}\\&=\dfrac{|\Hom(L'^\bullet,M^\bullet\oplus X^\bullet)^{Y^\bullet,\delta'_1,Z^\bullet[1],\delta'_2}_{L^\bullet,\delta'}|}{\langle \widehat{L'^\bullet[1]}, \widehat{L'_0}-\widehat{M_0}-\widehat{X_0}+\delta'\rangle}\cdot \dfrac{\{L'^\bullet[1],L'^\bullet\}  \{L'^\bullet[1],L^\bullet\} }{|\Aut(L'^\bullet)|\{L'^\bullet[1],M^\bullet\oplus X^\bullet\}}.
		\end{align*}
	\end{corollary}
	
	\begin{proposition}
		\label{prop:bijections-asso}
		Keep the notations as above. 
		There exist bijections
		\begin{align*}
			&\Hom(X^\bullet,L^\bullet)_{Y^\bullet,\delta_1}\times \Hom(M^\bullet,L^\bullet)_{Z^\bullet[1],\delta_2}\longrightarrow 
			\bigcup_{[L'],\delta}\Hom(M^\bullet\oplus X^\bullet,L^\bullet)^{Y^\bullet,\delta_1,Z^\bullet[1],\delta_2}_{L'^\bullet[1],\delta},
			\\ &\Hom(L'^\bullet,X^\bullet)_{Z^\bullet[1],\delta'_2}\times\Hom(L'^\bullet,M^\bullet)_{Y^\bullet,\delta'_1}
			\longrightarrow		\bigcup_{[L],\delta'}\Hom(L'^\bullet,M^\bullet\oplus X^\bullet)^{Y^\bullet,\delta'_1,Z^\bullet[1],\delta'_2}_{L^\bullet,\delta'}.
		\end{align*}
	\end{proposition}
	
	\begin{proof}
		The proof is similar to \cite[Proposition 3.5]{XX08}. We have a natural isomorphism
		\[
		\Hom(X^\bullet,L^\bullet)\times \Hom(M^\bullet,L^\bullet)\longrightarrow 
		\Hom(M^\bullet\oplus X^\bullet,L^\bullet).
		\]
		The first bijection follows by constraining on both sides of the map with the condition $$ 
		\cone f\cong Y, \widehat{\coker  H^0(f)}=\delta_1,
		\cone m\cong Z^\bullet[1], \widehat{\coker  H^0(m)}=\delta_2$$for $f:X^\bullet\rightarrow L^\bullet$ and $m:M^\bullet\rightarrow L^\bullet$. The second bijection is similar.
	\end{proof}
	
	Now we obtain the main result of this section.
	
	\begin{theorem}
		Let $\ca$ be a hereditary abelian $\bfk$-category with Euler form skew symmetric. Then the $\mathbb{Q}$-space $\och(\cR(\ca))$ is an associative algebra with $\{u_{[X^\bullet]}\mid X^\bullet\in\mathcal{R}(\ca)\}$ being the $\mathbb{Q}$-basis, and  the multiplication defined by
		\[
		u_{[Y^\bullet]}\cdot u_{[X^\bullet]}=\sum_{[L^\bullet]}F_{Y^\bullet X^\bullet}^{L^\bullet}\cdot u_{[L^\bullet]},
		\]
		and $u_{[0]}$ the identity.
	\end{theorem}
	
	\begin{proof}
		In order to prove the associativity, we need to prove 
		\begin{align}
			\label{eq:asso}
			\sum_{[L^\bullet]}F^{L^\bullet}_{Y^\bullet X^\bullet}F^{M^\bullet}_{L^\bullet Z^\bullet}=	\sum_{[L'^\bullet]}F^{L'^\bullet}_{X^\bullet Z^\bullet}F^{M^\bullet}_{Y^\bullet L'^\bullet},
		\end{align} 
		for any $X^\bullet, Y^\bullet, Z^\bullet$ and $M^\bullet$ in $\mathcal{R}(\ca)$.
		
		For simplicity, we denote $X^\bullet=X_0\oplus X_1[1]$ for $X_0,X_1\in\ca$, and similarly for others.
		
		Using Proposition \ref{prop:bijections-asso}, the left-hand side of \eqref{eq:asso} is equal to
		\begin{align*}
			&\sum_{[L^\bullet]}
			{\langle \widehat{Z_0},\widehat{L^\bullet}\rangle\langle \widehat{L_0},\widehat{Z^\bullet}\rangle }\sum_{\delta_2}\langle-\delta_2,\widehat{M^\bullet} \rangle
			\dfrac{|\Hom(M^\bullet,L^\bullet)_{Z^\bullet[1],\delta_2}|}{\langle \widehat{Z_0}+\widehat{L_0}-\widehat{M_0}-\delta_2,\widehat{L^\bullet}\rangle}\cdot \dfrac{\{L^\bullet[1],L^\bullet\}   }{|\Aut(L^\bullet)|\{M^\bullet[1],L^\bullet\}}\\\\&\quad\cdot
			{\langle \widehat{X_0},\widehat{Y^\bullet}\rangle\langle \widehat{Y_0},\widehat{X^\bullet}\rangle }\sum_{\delta_1}\langle -\widehat{Y_0}+\delta_1,\widehat{L^\bullet} \rangle
			\dfrac{|\Hom(X^\bullet,L^\bullet)_{Y^\bullet,\delta_1}|}{\langle \widehat{X^\bullet[1]}, \widehat{X_0}-\widehat{L_0}+\delta_1\rangle}\cdot \dfrac{\{X^\bullet[1],X^\bullet\}   }{|\Aut(X^\bullet)|\{X^\bullet[1],L^\bullet\}}
			\\
			&=
			\sum_{[L^\bullet],\delta_1,\delta_2}
			{\langle \widehat{Z_0},\widehat{L^\bullet}\rangle\langle \widehat{L_0},\widehat{Z^\bullet}\rangle\langle \widehat{X_0},\widehat{Y^\bullet}\rangle\langle \widehat{Y_0},\widehat{X^\bullet}\rangle }\dfrac{\langle-\delta_2,\widehat{M^\bullet} \rangle\langle -\widehat{Y_0}+\delta_1,\widehat{L^\bullet} \rangle}{\langle \widehat{Z_0}+\widehat{L_0}-\widehat{M_0}-\delta_2,\widehat{L^\bullet}\rangle\langle \widehat{X^\bullet[1]}, \widehat{X_0}-\widehat{L_0}+\delta_1\rangle}\\&\quad\cdot
			\dfrac{\{L^\bullet[1],L^\bullet\}   }{|\Aut(L^\bullet)|\{M^\bullet[1]\oplus X^\bullet[1],L^\bullet\}}	\dfrac{\{X^\bullet[1],X^\bullet\}   }{|\Aut(X^\bullet)|}\sum_{[L'^\bullet],\delta}|\Hom(M^\bullet\oplus X^\bullet,L^\bullet)^{Y^\bullet,\delta_1,Z^\bullet[1],\delta_2}_{L'^\bullet[1],\delta}|.	
		\end{align*}
		Dually, the right-hand side of \eqref{eq:asso} is equal to
		\begin{align*}
			&\sum_{[L'^\bullet]}
			{\langle \widehat{Z_0},\widehat{X^\bullet}\rangle\langle \widehat{X_0},\widehat{Z^\bullet}\rangle }\sum_{\delta'_2}\langle-\delta'_2,\widehat{L'^\bullet} \rangle
			\dfrac{|\Hom(L'^\bullet,X^\bullet)_{Z^\bullet[1],\delta'_2}|}{\langle \widehat{Z_0}+\widehat{X_0}-\widehat{L'_0}-\delta'_2,\widehat{X^\bullet}\rangle}\cdot \dfrac{\{X^\bullet[1],X^\bullet\}   }{|\Aut(X^\bullet)|\{L'^\bullet[1],X^\bullet\}}
			\\&\quad\cdot
			{\langle \widehat{L'_0},\widehat{Y^\bullet}\rangle\langle \widehat{Y_0},\widehat{L'^\bullet}\rangle }\sum_{\delta'_1}\langle -\widehat{Y_0}+\delta'_1,\widehat{M^\bullet} \rangle
			\dfrac{|\Hom(L'^\bullet,M^\bullet)_{Y^\bullet,\delta'_1}|}{\langle \widehat{L'^\bullet[1]}, \widehat{L'_0}-\widehat{M_0}+\delta'_1\rangle}\cdot \dfrac{\{L'^\bullet[1],L'^\bullet\}   }{|\Aut(L'^\bullet)|\{L'^\bullet[1],M^\bullet\}}
			\\&=
			\sum_{[L'^\bullet],\delta'_1,\delta'_2}
			{\langle \widehat{Z_0},\widehat{X^\bullet}\rangle\langle \widehat{X_0},\widehat{Z^\bullet}\rangle\langle \widehat{L'_0},\widehat{Y^\bullet}\rangle\langle \widehat{Y_0},\widehat{L'^\bullet}\rangle  }\dfrac{\langle-\delta'_2,\widehat{L'^\bullet} \rangle\langle -\widehat{Y_0}+\delta'_1,\widehat{M^\bullet} \rangle}{\langle \widehat{Z_0}+\widehat{X_0}-\widehat{L'_0}-\delta'_2,\widehat{X^\bullet}\rangle\langle \widehat{L'^\bullet[1]}, \widehat{L'_0}-\widehat{M_0}+\delta'_1\rangle}
			\\&\quad\cdot\dfrac{\{L'^\bullet[1],L'^\bullet\} }{|\Aut(L'^\bullet)|\{L'^\bullet[1],M^\bullet\oplus X^\bullet\}}	\dfrac{\{X^\bullet[1],X^\bullet\}   }{|\Aut(X^\bullet)|}\sum_{[L^\bullet],\delta'}|\Hom(L'^\bullet,M^\bullet\oplus X^\bullet)^{Y^\bullet,\delta'_1,Z^\bullet[1],\delta'_2}_{L^\bullet,\delta'}|.
		\end{align*}
		
		In order to prove the above two expressions are equal to each other, we relate variables such that {$\delta'=\widehat{L_0}-\delta$, $\delta'_1=\delta_1-\delta$ and $\delta'_2=\delta_2-\delta$.} By using Corollary  \ref{cor2}, it is enough to check the following equality:
		
		\begin{align}
			\label{3.14}
			\langle \widehat{L'_0}&-\widehat{X_0}-\delta,\widehat{L^\bullet}\rangle{\langle \widehat{L_0},\widehat{Z^\bullet}\rangle }{\langle \widehat{X_0},\widehat{Y^\bullet}\rangle\langle \widehat{Y_0},\widehat{X^\bullet}\rangle }\dfrac{\langle-\delta_2,\widehat{M^\bullet} \rangle\langle -\widehat{Y_0}+\delta_1,\widehat{L^\bullet} \rangle}{\langle -\delta_2,\widehat{L^\bullet}\rangle\langle \widehat{X_0}-\widehat{L_0}+\delta_1,\widehat{X^\bullet}\rangle}
			\\\notag
			&=\langle-\widehat{X_0}+\delta',\widehat{L'^\bullet}\rangle
			{\langle \widehat{X_0},\widehat{Z^\bullet}\rangle }{\langle \widehat{L'_0},\widehat{Y^\bullet}\rangle\langle \widehat{Y_0},\widehat{L'^\bullet}\rangle  }\dfrac{\langle-\delta'_2,\widehat{L'^\bullet} \rangle\langle -\widehat{Y_0}+\delta'_1,\widehat{M^\bullet} \rangle}{\langle \widehat{X_0}-\widehat{L'_0}-\delta'_2,\widehat{X^\bullet}\rangle\langle  \delta'_1,\widehat{L'^\bullet}\rangle},
		\end{align}
		by noting that $\widehat{X^\bullet[1]}=-\widehat{X^\bullet}$ and the Euler form of $\ca$ is skew symmetric.
		
		In fact, by replacing $\delta'$, $\delta_1'$, $\delta_2'$ with $\widehat{L_0}-\delta$, $\delta_1-\delta$ and $\delta_2-\delta$ in \eqref{3.14} respectively, we have
		\begin{align}
			\label{eq:Euler-simplify}
			\dfrac{\text{LHS\eqref{3.14}}}{\text{RHS}\eqref{3.14}}= &\langle \widehat{L_0}, \widehat{Z^\bullet}+\widehat{X^\bullet}-\widehat{L'^\bullet}\rangle
			\langle \widehat{X_0},\widehat{Y^\bullet}-\widehat{Z^\bullet}+\widehat{L'^\bullet}-\widehat{L^\bullet}\rangle\langle \widehat{L'_0},\widehat{L^\bullet} -\widehat{Y^\bullet}-\widehat{X^\bullet}\rangle
			\\
			\notag
			&\cdot\langle \widehat{Y_0}-\delta_1-\delta_2+\delta, \widehat{X^\bullet}-\widehat{L^\bullet}-\widehat{L'^\bullet}+\widehat{M^\bullet} \rangle.
		\end{align}
		Using the triangles in \eqref{eq:octahedral}, we have 
		\begin{align*}
			&\widehat{Z^\bullet}+\widehat{X^\bullet}-\widehat{L'^\bullet}=0,
			&\widehat{Y^\bullet}-\widehat{L^\bullet}=-\widehat{X^\bullet}= \widehat{Z^\bullet}-\widehat{L'^\bullet},
			\\
			&
			\widehat{L^\bullet }-\widehat{Y^\bullet}-\widehat{X^\bullet}=0,
			&\widehat{X^\bullet}-\widehat{L^\bullet}=-\widehat{Y^\bullet}=\widehat{L'^{\bullet}}-\widehat{M^\bullet},
		\end{align*}
		and then LHS\eqref{eq:Euler-simplify}$=1$.
		Then the desired formula \eqref{3.14} follows, which yields \eqref{eq:asso}.
		
		The proof is completed.
	\end{proof}

	\begin{remark}
		\label{rem}
		Inspired by the extended Hall algebra and also its Drinfeld double \cite{Rin90,Gr95,X97}, we can append a quantum torus to the derived Hall algebra of a root category, cf. \cite{Zh22}. Then this extended derived Hall algebra is well defined even the assumption of the Euler form of  $\ca$ being skew symmetric is dropped. 
		
		Moreover, the extended derived Hall algebra can be used to realize the quantum groups. But this construction is a little artificial. In fact, Bridgeland's Hall algebras \cite{Br13} or semi-derived Hall algebras \cite{LP21} are a better choice to realize quantum groups.
	\end{remark}

	\subsection{Twisted derived Hall algebras}
	
	For any $\delta\in K_0(\ca)$, $L^\bullet\cong L_0\oplus L_1[1]$, $M^\bullet\cong M_0\oplus M_1[1]$, $Z^\bullet\cong Z_0\oplus Z_1[1]$ in $\cR(\ca)$ with $L_i,M_i,Z_i\in\ca$ for $i=0,1$,  and considering triangles $M^\bullet\rightarrow L^\bullet\rightarrow Z^\bullet[1]\rightarrow M^\bullet[1]$, by Proposition \ref{prop4.1} and Lemma \ref{lem:W=W'}, we have 
	\begin{align}
		\label{3.17}
		&\quad	\dfrac{|\Hom(L^\bullet,Z^\bullet[1])_{M^\bullet[1],\delta''}|}{\langle \widehat{M_0}+\widehat{Z_1}-\widehat{L_0}-\delta'',\widehat{Z^\bullet[1]}\rangle}\dfrac{\{Z^\bullet,Z^\bullet[1]\}   }{|\Aut(Z^\bullet)|\{L^\bullet[1],Z^\bullet[1]\}}
		\\\notag
		&=	\dfrac{|\Hom(M^\bullet,L^\bullet)_{Z^\bullet[1],\delta}|}{\langle \widehat{M^\bullet[1]}, \widehat{M_0}-\widehat{L_0}+\delta\rangle}\dfrac{\{M^\bullet[1],M^\bullet\}   }{|\Aut(M^\bullet)|\{M^\bullet[1],L^\bullet\}},
	\end{align}where $\delta''=\widehat{Z_1}-\delta$. Substituting (\ref{3.17}) into  (\ref{3.10}), we get the following formula
	\begin{align}
		\label{3.18}
		F_{L^\bullet,Z^\bullet}^{M^\bullet}
		&={\langle \widehat{Z_0},\widehat{L^\bullet}\rangle\langle \widehat{L_0},\widehat{Z^\bullet}\rangle }\sum_{\delta}\langle-\delta,\widehat{M^\bullet} \rangle
		\dfrac{|\Hom(L^\bullet,Z^\bullet[1])_{M^\bullet[1],\delta''}|}{\langle \widehat{Z_0}+\widehat{L_0}-\widehat{M_0}-\delta,\widehat{L^\bullet}\rangle}\cdot\dfrac{\{L^\bullet[1],L^\bullet\}   }{\{M^\bullet[1],L^\bullet\}}
		\\\notag &\quad\cdot 
		\dfrac{\{Z^\bullet,Z^\bullet[1]\} |\Aut(M^\bullet)|\{M^\bullet[1],L^\bullet\}  }{|\Aut(L^\bullet)||\Aut(Z^\bullet)|\{L^\bullet[1],Z^\bullet[1]\}\{M^\bullet[1],M^\bullet\}}\cdot
		\dfrac{\langle \widehat{M^\bullet[1]}, \widehat{M_0}-\widehat{L_0}+\delta\rangle}{\langle \widehat{M_0}+\widehat{Z_1}-\widehat{L_0}-\delta'',\widehat{Z^\bullet[1]}\rangle}
		\\\notag&=
		\dfrac{\{L^\bullet[1],L^\bullet\} \{Z^\bullet,Z^\bullet[1]\} |\Aut(M^\bullet)|   }{|\Aut(L^\bullet)||\Aut(Z^\bullet)|\{M^\bullet[1],M^\bullet\}}
		\sum_{\delta}\dfrac{|\Hom(L^\bullet,Z^\bullet[1])_{M^\bullet[1],\widehat{Z_1}-\delta}|}{\{L^\bullet[1],Z^\bullet[1]\}}
		\\ \notag&\quad\cdot 
		\dfrac{	\langle-\delta,\widehat{M^\bullet} \rangle{\langle \widehat{Z_0},\widehat{L^\bullet}\rangle\langle \widehat{L_0},\widehat{Z^\bullet}\rangle }\langle \widehat{M^\bullet[1]}, \widehat{M_0}-\widehat{L_0}+\delta\rangle}{\langle \widehat{Z_0}+\widehat{L_0}-\widehat{M_0}-\delta,\widehat{L^\bullet}\rangle\langle \widehat{M_0}-\widehat{L_0}+\delta,\widehat{Z^\bullet[1]}\rangle}
		\\\notag&=
		\dfrac{\{L^\bullet[1],L^\bullet\} \{Z^\bullet,Z^\bullet[1]\} |\Aut(M^\bullet)|   }{|\Aut(L^\bullet)||\Aut(Z^\bullet)|\{M^\bullet[1],M^\bullet\}}
		\sum_{\delta}\dfrac{	\langle\delta,\widehat{M^\bullet} \rangle\langle  \widehat{M_0},2\widehat{M^\bullet}\rangle}{\langle \widehat{L_0},\widehat{M^\bullet}+\widehat{L^\bullet}\rangle}
		\\\notag&\quad\cdot
		\dfrac{|\Hom(L^\bullet,Z^\bullet[1])_{M^\bullet[1],\widehat{Z_1}-\delta}|}{\{L^\bullet[1],Z^\bullet[1]\}}.
	\end{align}
	Observe that \[
	\dfrac{\{Z^\bullet,Z^\bullet[1]\}}{\{Z^\bullet[1],Z^\bullet\}}
	=\dfrac{|\Hom_{\ca}(Z_0,Z_1)|\cdot|\Ext^1_{\ca}(Z_1,Z_0)|}{|\Hom_{\ca}(Z_1,Z_0)|\cdot|\Ext^1_{\ca}(Z_0,Z_1)|}
	=\dfrac{\langle \widehat{Z_0},\widehat{Z_1}\rangle}{\langle \widehat{Z_1},\widehat{Z_0}\rangle}
	=\langle \widehat{Z_0},2\widehat{Z_1}\rangle,
	\]and
	\[
	\{L^\bullet[1],Z^\bullet[1]\}=|\Hom_{\ca}(L_1,Z_1)|\cdot|\Ext^1_{\ca}(L_0,Z_0)|=\dfrac{|\Hom_{\ca}(L_1,Z_1)|\cdot|\Hom_{\ca}(L_0,Z_0)|}{\langle L_0,Z_0\rangle}.
	\]
	Hence formula (\ref{3.18}) is equal to
	\begin{align}
		\label{eq:Hallmulti-reform2}
		&	\dfrac{\{L^\bullet[1],L^\bullet\} \{Z^\bullet[1],Z^\bullet\} |\Aut(M^\bullet)|   }{|\Aut(L^\bullet)||\Aut(Z^\bullet)|\{M^\bullet[1],M^\bullet\}} \cdot\langle \widehat{Z_0},2\widehat{Z_1}\rangle\langle\widehat{L_0},\widehat{Z_0}\rangle
		\\\notag &\qquad\cdot\sum_{\delta}\dfrac{	\langle\delta,\widehat{M^\bullet} \rangle\langle  \widehat{M_0},2\widehat{M^\bullet}\rangle}{\langle \widehat{L_0},\widehat{M^\bullet}+\widehat{L^\bullet}\rangle}
		\dfrac{|\Hom(L^\bullet,Z^\bullet[1])_{M^\bullet[1],\widehat{Z_1}-\delta}|}{|\Hom(L_1,Z_1)||\Hom(L_0,Z_0)|}.
	\end{align}
	
	We define $\ch(\cR(\ca))$ to be the same vector space as  $\och(\cR(\ca))$, with the multiplication defined by
	\begin{align}
		\label{3.20}
		u_{[L^\bullet]}\diamond u_{[Z^\bullet]}=\sum_{[M^\bullet]}\langle \widehat{Z_0}&,2\widehat{Z_1}\rangle\langle\widehat{L_0},\widehat{Z_0}\rangle \sum_{\delta}\dfrac{	\langle\delta,\widehat{M^\bullet} \rangle\langle  \widehat{M_0},2\widehat{M^\bullet}\rangle }{\langle \widehat{L_0},\widehat{M^\bullet}+\widehat{L^\bullet}\rangle}
		\\\notag&\quad\cdot
		\dfrac{|\Hom(L^\bullet,Z^\bullet[1])_{M^\bullet[1],\widehat{Z_1}-\delta}|}{|\Hom(L_1,Z_1)||\Hom(L_0,Z_0)|}\cdot u_{[M^\bullet]}.
	\end{align}
	It follows from \eqref{eq:Hallmulti-reform2} that $\ch(\cR(\ca))$ is an associative algebra, and there is an isomorphism $\och(\cR(\ca))\stackrel{\cong}{\longrightarrow}\ch(\cR(\ca))$
	by sending $u_{[L^\bullet]}\mapsto u_{[L^\bullet]}\cdot \dfrac{\{L^\bullet[1],L^\bullet\}}{|\Aut(L^\bullet)|}$.
	
	The algebra $\ch(\cR(\ca))$ is called the Drinfeld dual of $\och(\cR(\ca))$.

	Similar to the Ringel-Hall algebras \cite{Rin90,Gr95}, we construct a twisted version of the derived Hall algebra $\ch(\cR(\ca))$. 
	Note that $\widehat{Z^\bullet}+\widehat{L^\bullet}=\widehat{M^\bullet}$ holds for any triangle $Z^\bullet\rightarrow M^\bullet\rightarrow L^\bullet\rightarrow Z^\bullet[1]$. Let $\sqq=\sqrt{q}$. Hence we define the twisted derived Hall algebra  $\widetilde{\ch}(\cR(\ca))$ to be the  $\Q(\sqq)$-linear space with the same basis as the Drinfeld dual derived Hall algebra, and the multiplication given by
	\begin{align}
		u_{[L^\bullet]}*u_{[Z^\bullet]}={\langle \widehat{L^\bullet},\widehat{Z^\bullet}\rangle^{1/2}}\cdot u_{[L^\bullet]}\diamond u_{[Z^\bullet]}.
	\end{align}

	\section{Drinfeld double of Hall algebras}
	\label{sec:Drinfelddouble}
	
	\subsection{Hall algebras}

	Recall that $\ca$ is a hereditary abelian $\bfk$-linear category which is essentially small with finite-dimensional morphism and extension spaces. We also assume the Euler form of $\ca$ to be skew symmetric.
	
	Given objects $M,N,L\in\ca$, define $\Ext^1(M,N)_L\subseteq \Ext^1(M,N)$ as the subset parameterizing extensions whose middle term is isomorphic to $L$. We define the Hall algebra $\ch(\ca)$  to be the $\Q$-linear space whose basis is formed by the isoclasses $[M]$ of objects $M$ in $\ca$, with the multiplication defined by (see \cite{Br13})
	\begin{align}
		\label{eq:mult}
		[M]\diamond [N]=\sum_{[L]\in \Iso(\ca)}\frac{|\Ext^1(M,N)_L|}{|\Hom(M,N)|}[L].
	\end{align}
	
	For any three objects $L,M,N$, let
	\begin{align}
		\label{eq:Fxyz}
		G_{MN}^L:= \big |\{X\subseteq L \mid X \cong N,  L/X\cong M\} \big |.
	\end{align}
	The Riedtmann-Peng formula states that
	\[
	G_{MN}^L= \frac{|\Ext^1(M,N)_L|}{|\Hom(M,N)|} \cdot \frac{|\Aut(L)|}{|\Aut(M)| |\Aut(N)|}.
	\]
	The Hall numbers satisfy the following associativity for any $L,M,N,Z\in \mathcal{A}$:
	\begin{align*}
		\sum_{[X]\in{\rm{Iso}}(\mathcal{A})} G_{LM}^{X}G_{XN}^{Z}
		=\sum_{[X']\in{\rm{Iso}}(\mathcal{A})} G_{LX'}^{Z}G_{MN}^{X'}:=G_{LMN}^Z.
	\end{align*}

	
	The twisted Hall algebra $\widetilde{\ch}(\ca)$ is defined over $\ch(\ca)\otimes_\Q \Q(\sqq)$ with multiplication twisted by Euler form:
	\begin{align*}
		[M]*[N]={\langle M,N \rangle^{1/2}}\cdot [M]\diamond[N].
	\end{align*}
	
	We denote by $\widetilde{\ch}(\ca)\widehat{\otimes} \widetilde{\ch}(\ca)$ the space of formal linear combinations
	$$\sum_{[A],[B]\in \Iso(\ca)} c_{A,B} [A]\otimes[B],$$
	where $\widehat{\otimes}$ is the {completed tensor product}.
	The coproduct and counit for $\widetilde{\ch}(\ca)$ are given by Green \cite{Gr95} (see also \cite{X97}):
	\begin{align*}
		\Delta:&\widetilde{\ch}(\ca)\longrightarrow \widetilde{\ch}(\ca)\widehat{\otimes} \widetilde{\ch}(\ca),\qquad \epsilon: \widetilde{\ch}(\ca)\longrightarrow \Q(\sqq),
	\end{align*}
	such that
	\begin{align*}
		\Delta([A])=&\sum_{[B],[C]}\langle \widehat{B},\widehat{C}\rangle^{1/2}G_{B,C}^A [B]\otimes[C],
		\qquad		\epsilon([A])=\delta_{[A],0},
	\end{align*}
	for $[A]\in\Iso(\ca)$. Then $(\widetilde{\ch}(\ca),*,[0],\Delta,\epsilon)$ is a \emph{topological bialgebra} defined over $\Q(\sqq)$ (see \cite{Gr95,X97}).
	Here ``topological'' means that everything should be considered in the completed space.
	\begin{remark}
		If $\ca$ is in particular a finite length hereditary category (for instance, the category of nilpotent finite-dimensional representations of a finite quiver), then $(\widetilde{\ch}(\ca),*,[0],\Delta,\epsilon)$ is a genuine bialgebra. Moreover, Xiao \cite{X97} defined the antipode which endowed $\widetilde{\ch}(\ca)$ with a natural Hopf algebra structure.
	\end{remark}

	\subsection{Drinfeld double}
	The bilinear pairing $(-,-):\widetilde{\ch}(\ca)\times \widetilde{\ch}(\ca)\rightarrow\Q(\sqq)$ given by
	\begin{equation}
		([M],[N])=\delta_{[M],[N]}|\Aut(M)|
	\end{equation}
	is a {Hopf pairing} on $\widetilde{\ch}(\ca)$; see \cite{Gr95}.


	For the hereditary abelian category $\ca$, there is a unique algebra structure on $\widetilde{\ch}(\ca){\otimes}\widetilde{\ch}(\ca)$ satisfying the following conditions, which is called the {Drinfeld double} of $\widetilde{\ch}(\ca)$, and denoted by $\cd\th(\ca)$; see \cite{X97,Cr10,LP21}.
	
	\begin{itemize}
		\item[(D1)] The maps $$\widetilde{\ch}(\ca)\longrightarrow\cd\th(\ca),\quad a\mapsto a\otimes1$$
		and
		$$\widetilde{\ch}(\ca)\longrightarrow\cd\th(\ca),\quad a\mapsto 1\otimes a$$
		are injective homomorphisms of $\Q(\sqq)$-algebras. 
		\item[(D2)] For any $a,b\in \widetilde{\ch}(\ca)$, one has
		$$(a\otimes1)(1\otimes b)=a\otimes b.$$
		\item[(D3)] For any $a,b\in \widetilde{\ch}(\ca)$, one has
		\begin{equation}\label{equation drinfeld double 2}
			\sum (a_{(2)},b_{(1)})a_{(1)}\otimes b_{(2)}=\sum (a_{(1)},b_{(2)}) (1\otimes b_{(1)})(a_{(2)}\otimes1).
		\end{equation}
		Here we use Sweedler’s notation: $\Delta(a)=\sum a_{(1)}\otimes a_{(2)}$.
	\end{itemize}
	
	It is worth noting that if $\widetilde{\ch}(\ca)$ is a topological bialgebra, then one should replace the tensor product $\otimes$
	in the statement by the completed one
	$\widehat{\otimes}$.
	
	In particular, if $\ca$ is a finite length hereditary abelian category, then the Drinfeld double of $\widetilde{\ch}(\ca)$ is again a Hopf algebra.

	\subsection{Realization of Drinfeld double of Hall algebras}
	
	First, we shall give a basis of the twisted derived Hall algebra $\widetilde{\ch}(\cR(\ca))$ of the root category $\cR(\ca)$. 
	
	\begin{lemma}
		\label{basis}
		The algebra $\widetilde{\ch}(\cR(\ca))$ has a basis given by
		\begin{align}
			\label{eq:monomialbasis}
			\{u_{[A]}*u_{[B[1]]}\mid A,B\in\ca\}.
		\end{align}
	\end{lemma}
	
	\begin{proof}
		This proof is inspired by \cite[Theorem 2.25]{LP21}; see also \cite[Lemma 2.1]{GP} and \cite[Proposition 4.8]{BG16}.
		
		From the definition, $\widetilde{\ch}(\cR(\ca))$ has a basis $\{u_{[A\oplus B[1]]}\mid A,B\in\ca\}$. We have the formula
		\begin{align*}
			u_{[A]}*u_{[B[1]]}&=\sum_{[M^\bullet]}\langle \widehat{A},\widehat{B}\rangle^{-1/2} \dfrac{	\langle\widehat{A}-\widehat{M_0},\widehat{M^\bullet} \rangle\langle  \widehat{M_0},2\widehat{M^\bullet}\rangle }{\langle \widehat{A},\widehat{M^\bullet}+\widehat{A}\rangle}\cdot |\Hom(A,B)_{M^\bullet[1]}|\cdot u_{[M^\bullet]}
			\\\notag&=\sum_{[M^\bullet]} \dfrac{1}{\langle \widehat{A},\widehat{B}\rangle^{1/2}\langle \widehat{M_0},\widehat{M_1}\rangle}\cdot |\Hom(A,B)_{M^\bullet[1]}|\cdot u_{[M^\bullet]},
		\end{align*}
		where $M^\bullet\cong M_0\oplus M_1[1]$ in $\cR(\ca)$ with $M_i\in\ca$ for $i=0,1$, such that there exists a triangle in $\cR(\ca)$
		\begin{equation}
			\label{eq:trianglebasis}
			B[1]\stackrel{l}{\longrightarrow} M^\bullet \stackrel{m}{\longrightarrow} A\stackrel{n}{\longrightarrow} B.
		\end{equation}
		
		Applying $(M^\bullet,-):=\Hom_{\cR(\ca)}(M^\bullet,-)$ to \eqref{eq:trianglebasis}, we have
		\[\cdots \longrightarrow (M^\bullet,B[1]) \longrightarrow \End_{\cR(\ca)} (M^\bullet) \longrightarrow (M^\bullet, A)\longrightarrow\cdots.\]
		Hence $|\End_{\cR{(\ca)}}(M^\bullet)|\le|(M^\bullet,B[1])| |(M^\bullet,A)|=|(M^\bullet,A\oplus B[1])|$. Similarly, applying $(-,A\oplus B[1])$, we have a long exact sequence
		\[
		\begin{tikzpicture}
			\node (0) at (-5*1.3,0) {$\cdots$};
			\node (-2) at (-3*1.3,0) {$(A,A\oplus B[1])$};
			\node (2) at (0,0) {$(M^\bullet,A\oplus B[1])$};
			\draw[->] (-2) --node[above ]{} (2);
			\draw[->] (0) --node[above ]{} (-2);
			\node (3) at (3*1.3,0) {$(B[1],A\oplus B[1])$};
			\draw[->] (2) --node[above ]{} (3);
			\node (4) at (6.5,0) {$\cdots$,};
			\draw[->] (3) --node[above ]{} (4);
		\end{tikzpicture}
		\]
		and then$|(M^\bullet,A\oplus B[1])|\le|\End_{\cR{(\ca)}}(A\oplus B[1])| $. Therefore, we get
		$$|\End_{\cR{(\ca)}}(A\oplus B[1])|\ge |\End_{\cR{(\ca)}}(M^\bullet)|= |\End_{\cR{(\ca)}}(M_0\oplus M_1[1])|.$$ In addition, if the equality holds, we have that 
		$$(n[1],A\oplus B[1]):(B[1],A\oplus B[1])\rightarrow (A[1],A\oplus B[1])=0.$$ 
		It follows that $n[1]=0$ which implies that $A\oplus B[1]\cong M^\bullet\cong M_0\oplus M_1[1]$.
		Using the formula
		\begin{align*}
			\dfrac{ u_{[A\oplus B[1]]}}{\langle \widehat{A},\widehat{B}\rangle^{3/2}}    &=u_{[A]}*u_{[B[1]]}-\sum_{[M^\bullet],M^\bullet\ncong A\oplus B[1]} \dfrac{1}{\langle \widehat{A},\widehat{B}\rangle^{1/2}\langle \widehat{M_0},\widehat{M_1}\rangle}\cdot |\Hom(A,B)_{M^\bullet[1]}|\cdot u_{[M^\bullet]},  
		\end{align*}
		and $|\End_{\cR{(\ca)}}(A\oplus B[1])|<\infty$, we can obtain by induction that $u_{[A\oplus B[1]]}$ can be written as a linear combination of the set $\{u_{[A]}*u_{[B[1]]}\mid A,B\in\ca\}$.
		
		On the other hand, it is not hard to prove that $\{u_{[A]}*u_{[B[1]]}\mid A,B\in\ca\}$ is linearly independent. Therefore, the set $\{u_{[A]}*u_{[B[1]]}\mid A,B\in\ca\}$ is a basis of $\widetilde{\ch}(\cR(\ca))$.
	\end{proof}

	In the following, we prove that the twisted derived Hall algebra $\th(\cR(\ca))$ is isomorphic to the Drinfeld double of $\widetilde{\ch}(\ca)$. When $\widetilde{\ch}(\ca)$ is only a topological bialgebra, its Drinfeld double is the completed tensor product $\widetilde{\ch}(\ca)\widehat{\otimes} \widetilde{\ch}(\ca)$ (as a vector space). Accordingly, in this case, we need to complete the twisted derived Hall algebra with respect to the basis
	\eqref{eq:monomialbasis}.

	\begin{lemma}
		\label{embed}
		For a hereditary abelian $\bfk$-category $\ca$ with Euler form skew symmetric, there exist two embeddings of algebras
		\begin{align*}
			\Psi^+:\widetilde{\ch}&(\ca)\longrightarrow \th(\cR(\ca)),
			\\&[M]\mapsto u_{[M]};
		\end{align*}
		and
		\begin{align*}
			\Psi^-:\widetilde{\ch}&(\ca)\longrightarrow \th(\cR(\ca)),
			\\&[M]\mapsto u_{[M[1]]}.
		\end{align*}
	\end{lemma}
	
	\begin{proof}
		It is clear that $\Psi^+$ and $\Psi^-$ are injections by Lemma \ref{basis}. 
		
		For any $M^\bullet$, $L^\bullet\cong L_0$, $Z^\bullet\cong Z_0$ in $\cR(\ca)$ with $L_0,Z_0\in\ca$, if there is a triangle $Z^\bullet \xrightarrow{l}M^\bullet\xrightarrow{m} L^\bullet\xrightarrow{n}Z^\bullet[1]$, then $M^\bullet$ is isomorphic to $M_0$ with $M_0\in\ca$ and $\delta=\widehat{\coker  H^0(m)}$ is zero. Therefore, we have
		\begin{align}
			\label{3.22}
			u_{[L_0]}*u_{[Z_0]}&=\sum_{[M_0]}\langle \widehat{L_0},\widehat{Z_0}\rangle^{1/2}\cdot
			\dfrac{|\Hom(L_0,Z_0[1])_{M_0[1]}|}{|\Hom(L_0,Z_0)|}\cdot u_{[M_0]}
			\\\notag&=\sum_{[M_0]}\langle \widehat{L_0},\widehat{Z_0}\rangle^{1/2}\cdot
			\dfrac{|\Ext^1_\ca(L_0,Z_0)_{M_0}|}{|\Hom(L_0,Z_0)|}\cdot u_{[M_0]}.
		\end{align}
		Similarly, for any $L^\bullet\cong L_1[1]$, $Z^\bullet\cong Z_1[1]$ in $\cR(\ca)$ with $L_1,Z_1\in\ca$ we have
		\begin{align}
			\label{3.23}
			u_{[L_1[1]]}*u_{[Z_1[1]]}&=\sum_{[M_1[1]]}\langle \widehat{L_1},\widehat{Z_1}\rangle^{1/2}\cdot
			\dfrac{|\Hom(L_1[1],Z_1)_{M_1}|}{|\Hom(L_1,Z_1)|}\cdot u_{[M_1[1]]}
			\\\notag&=\sum_{[M_1[1]]}\langle \widehat{L_1},\widehat{Z_1}\rangle^{1/2}\cdot
			\dfrac{|\Ext^1_\ca(L_1,Z_1)_{M_1}|}{|\Hom(L_1,Z_1)|}\cdot u_{[M_1[1]]}.
		\end{align}
		
		It follows directly from the formulas (\ref{3.22}) and (\ref{3.23}) that $\Psi^+$ and $\Psi^-$ are embeddings.
	\end{proof}

	\begin{theorem}
		\label{iso}
		For a hereditary abelian $\bfk$-category $\ca$ with Euler form skew symmetric, the twisted derived Hall algebra $\widetilde{\ch}(\cR(\ca))$ is isomorphic to the Drinfeld double $\cd\th(\ca)$.
	\end{theorem}

	\begin{proof}
		By Lemma \ref{basis}, we have an isomorphism of $\mathbb{Q}(\sqq)$-vector spaces
		\begin{align*}
			\Psi_\ca:&\cd\widetilde{\ch}(\ca)\longrightarrow \widetilde{\ch}(\cR(\ca)),
			\\&[M_0]\otimes[M_1]\mapsto u_{[M_0]}*u_{[{M_1[1]}]}.
		\end{align*}
		Due to Lemma \ref{embed}, it remains to prove the relation
		\begin{align}
			\label{dd}
			\sum(x_{(2)},y_{(1)})\cdot\Psi^+(x_{(1)})*\Psi^-(y_{(2)})=\sum(x_{(1)},y_{(2)})\cdot\Psi^-(y_{(1)})*\Psi^+(x_{(2)})
		\end{align}
		for $x,y\in \widetilde{\ch}(\ca)$, 
		where $\Delta(x)=\sum x_{(1)}\otimes x_{(2)}$ and $\Delta(y)=\sum y_{(1)}\otimes y_{(2)}$. 
		
		We do some direct calculation before confirming the relation.
		For $M^\bullet\cong M_0\oplus M_1[1]$, $L^\bullet\cong L_0$, $Z^\bullet\cong Z_1[1]$ in $\cR(\ca)$ with $M_0,M_1,L_0,Z_1\in\ca$, and a triangle $Z^\bullet \xrightarrow{l}M^\bullet\xrightarrow{m} L^\bullet\xrightarrow{n}Z^\bullet[1]$, noting that $\delta=\widehat{\coker  H^0(m)}$ is equal to $\widehat{L_0}-\widehat{M_0}$, it follows that
		\begin{align}
			\label{3.24}
			u_{[L_0]}*u_{[Z_1[1]]}&=\sum_{[M^\bullet]}\langle \widehat{L_0},\widehat{Z_1}\rangle^{-1/2} \dfrac{	\langle\widehat{L_0}-\widehat{M_0},\widehat{M^\bullet} \rangle\langle  \widehat{M_0},2\widehat{M^\bullet}\rangle }{\langle \widehat{L_0},\widehat{M^\bullet}+\widehat{L^\bullet}\rangle}\cdot |\Hom(L_0,Z_1)_{M^\bullet[1]}|\cdot u_{[M^\bullet]}
			\\\notag&=\sum_{[M^\bullet]} \dfrac{1}{\langle \widehat{L_0},\widehat{Z_1}\rangle^{1/2}\langle \widehat{M_0},\widehat{M_1}\rangle}\cdot |\Hom(L_0,Z_1)_{M^\bullet[1]}|\cdot u_{[M^\bullet]}.
		\end{align}
		Similarly, for $L^\bullet\cong {L_1[1]}$, $Z^\bullet\cong {Z_0}$ in $\cR(\ca)$ with {$L_1,Z_0\in\ca$}, we also have
		\begin{align}
			\label{3.25}
			u_{[L_1[1]]}*u_{[Z_0]}&=\sum_{[M^\bullet]} \dfrac{1}{\langle \widehat{L_1},\widehat{Z_0}\rangle^{1/2}\langle \widehat{M_0},\widehat{M_1}\rangle^2}\cdot |\Hom(L_1[1],Z_0[1])_{M^\bullet[1]}|\cdot u_{[M^\bullet]}.
		\end{align}

		Now assume $x=[X_0]$ and $y=[X_1]$ with
		\[\Delta([X_0])=\sum_{[L_0],[Z_0]}\langle \widehat{L_0},\widehat{Z_0}\rangle^{1/2}G_{L_0Z_0}^{X_0}[L_0]\otimes [Z_0]\]
		and
		\[\Delta([X_1])=\sum_{[L_1],[Z_1]}\langle \widehat{L_1},\widehat{Z_1}\rangle^{1/2}G_{L_1Z_1}^{X_1}[L_1]\otimes [Z_1].\]
		Then we have \begin{align*}
			\text{LHS(\ref{dd})}&=\sum_{[L_0],[Z_0],[L_1],[Z_1]}\langle \widehat{L_0},\widehat{Z_0}\rangle^{1/2} \langle \widehat{L_1},\widehat{Z_1}\rangle^{1/2}G_{L_0Z_0}^{X_0} G_{L_1Z_1}^{X_1} ([Z_0],[L_1]) \cdot u_{[L_0]}*u_{[Z_1[1]]}
			\\&=
			\sum_{[L_0],[Z_0],[Z_1]}\langle \widehat{L_0},\widehat{Z_0}\rangle^{1/2} \langle \widehat{Z_0},\widehat{Z_1}\rangle^{1/2}G_{L_0Z_0}^{X_0} G_{Z_0Z_1}^{X_1} a_{Z_0} \cdot u_{[L_0]}*u_{[Z_1[1]]}
			\\&=
			\sum_{[L_0],[Z_0],[Z_1]}\langle \widehat{L_0},\widehat{Z_0}\rangle^{1/2} \langle \widehat{Z_0},\widehat{Z_1}\rangle^{1/2}G_{L_0Z_0}^{X_0} G_{Z_0Z_1}^{X_1} a_{Z_0} \cdot \sum_{[M^\bullet]} \dfrac{1}{\langle \widehat{L_0},\widehat{Z_1}\rangle^{1/2}\langle \widehat{M_0},\widehat{M_1}\rangle}
			\\&\quad\cdot |\Hom(L_0,Z_1)_{M^\bullet[1]}|\cdot u_{[M^\bullet]},
		\end{align*}
		where we use the formula (\ref{3.24}).
		
		It is well known that (see e.g. \cite[(8.8)]{SV})
		\[
		|\{l\in \Hom_\ca(L_0,Z_1)\mid \coker l\cong M_1, \ker l\cong M_0\}|=\sum_{[I]}G_{IM_0}^{L_0}G_{M_1I}^{Z_1}\cdot a_I.
		\]
		Hence, for $M^\bullet\cong M_0\oplus M_1[1]$ $\in\cR(\ca)$ with $M_0,M_1\in\ca$, we have
		\begin{align*}
			\text{LHS(\ref{dd})}&=\sum_{[L_0],[Z_0],[Z_1]}\langle \widehat{L_0},\widehat{Z_0}\rangle^{1/2} \langle \widehat{Z_0},\widehat{Z_1}\rangle^{1/2}G_{L_0Z_0}^{X_0} G_{Z_0Z_1}^{X_1} a_{Z_0} \cdot \sum_{[M^\bullet]} \dfrac{1}{\langle \widehat{L_0},\widehat{Z_1}\rangle^{1/2}\langle \widehat{M_0},\widehat{M_1}\rangle}
			\\&\quad\cdot \sum_{[I_0]}G_{I_0M_0}^{L_0}G_{M_1I_0}^{Z_1}\cdot a_{I_0}\cdot  u_{[M^\bullet]}
			\\&=
			\sum_{\substack{[L_0],[Z_0],[Z_1]\\ [M_0],[M_1],[I_0]}}(\dfrac{\langle \widehat{M_0}+\widehat{I_0},\widehat{Z_0}\rangle\langle \widehat{Z_0},\widehat{M_1}+\widehat{I_0}\rangle}{\langle {\widehat{M_0}+\widehat{I_0}},\widehat{M_1}+\widehat{I_0}\rangle})^{1/2}\dfrac{1}{\langle \widehat{M_0},\widehat{M_1}\rangle}
			\cdot a_{Z_0}a_{I_0}  
			\\&\quad\cdot G_{L_0Z_0}^{X_0} G_{Z_0Z_1}^{X_1}G_{I_0M_0}^{L_0}G_{M_1I_0}^{Z_1}\cdot  u_{[M^\bullet]}
			\\&=\sum_{\substack{[Z_0], [M_0],[M_1],[I_0]}}(\dfrac{\langle \widehat{M_0},\widehat{Z_0}\rangle\langle \widehat{Z_0},\widehat{M_1}\rangle}{\langle {\widehat{M_0}},\widehat{I_0}\rangle\langle {\widehat{I_0}},\widehat{M_1}\rangle})^{1/2}\dfrac{1}{\langle \widehat{M_0},\widehat{M_1}\rangle^{3/2}}
			\cdot a_{Z_0}a_{I_0}  
			\\&\quad\cdot G_{I_0M_0Z_0}^{X_0} G_{Z_0M_1I_0}^{X_1}\cdot  u_{[M^\bullet]}.
		\end{align*}
		
		Dually, we have
		\begin{align*}
			\text{RHS(\ref{dd})}&=\sum_{[L_0],[Z_0],[L_1],[Z_1]}\langle \widehat{L_0},\widehat{Z_0}\rangle^{1/2} \langle \widehat{L_1},\widehat{Z_1}\rangle^{1/2}G_{L_0Z_0}^{X_0} G_{L_1Z_1}^{X_1} ([L_0],[Z_1]) \cdot u_{[L_1[1]]}*u_{[Z_0]}
			\\&=
			\sum_{[L_0],[Z_0],[L_1]}\langle \widehat{L_0},\widehat{Z_0}\rangle^{1/2} \langle \widehat{L_1},\widehat{L_0}\rangle^{1/2}G_{L_0Z_0}^{X_0} G_{L_1L_0}^{X_1} a_{L_0} \cdot u_{[L_1[1]]}*u_{[Z_0]}
			\\&=
			\sum_{[L_0],[Z_0],[L_1]}\langle \widehat{L_0},\widehat{Z_0}\rangle^{1/2} \langle \widehat{L_1},\widehat{L_0}\rangle^{1/2}G_{L_0Z_0}^{X_0} G_{L_1L_0}^{X_1} a_{L_0} \cdot \sum_{[M^\bullet]} \dfrac{1}{\langle \widehat{L_1},\widehat{Z_0}\rangle^{1/2}\langle \widehat{M_0},\widehat{M_1}\rangle^2}
			\\&\quad\cdot |\Hom(L_1[1],Z_0[1]
			)_{M^\bullet[1]}|\cdot u_{[M^\bullet]},
		\end{align*}
		where we use the formula (\ref{3.25}).
		Hence
		\begin{align*}
			\text{RHS(\ref{dd})}&=\sum_{[L_0],[Z_0],[L_1]}\langle \widehat{L_0},\widehat{Z_0}\rangle^{1/2} \langle \widehat{L_1},\widehat{L_0}\rangle^{1/2}G_{L_0Z_0}^{X_0} G_{L_1L_0}^{X_1} a_{L_0} \cdot \sum_{[M^\bullet]} \dfrac{1}{\langle \widehat{L_1},\widehat{Z_0}\rangle^{1/2}\langle \widehat{M_0},\widehat{M_1}\rangle^2}
			\\&\quad\cdot \sum_{[I_1]}G_{I_1M_1}^{L_1}G_{M_0I_1}^{Z_0}\cdot a_{I_1}\cdot  u_{[M^\bullet]}
			\\&=
			\sum_{\substack{[L_0],[Z_0],[L_1]\\ [M_0],[M_1],[I_1]}}(\dfrac{\langle \widehat{L_0},\widehat{M_0}+\widehat{I_1}\rangle\langle \widehat{M_1}+\widehat{I_1},\widehat{L_0}\rangle}{\langle {\widehat{M_1}+\widehat{I_1}},\widehat{M_0}+\widehat{I_1}\rangle})^{1/2}\dfrac{1}{\langle \widehat{M_0},\widehat{M_1}\rangle^2}
			\cdot a_{L_0}a_{I_1}  
			\\&\quad\cdot G_{L_0Z_0}^{X_0} G_{L_1L_0}^{X_1}G_{I_1M_1}^{L_1}G_{M_0I_1}^{Z_0}\cdot  u_{[M^\bullet]}
			\\&=\sum_{\substack{[L_0], [M_0],[M_1],[I_1]}}(\dfrac{\langle \widehat{L_0},\widehat{M_0}\rangle\langle \widehat{M_1},\widehat{L_0}\rangle}{\langle {\widehat{M_1}},\widehat{I_1}\rangle\langle {\widehat{I_1}},\widehat{M_0}\rangle})^{1/2}\dfrac{1}{\langle \widehat{M_0},\widehat{M_1}\rangle^{3/2}}
			\cdot a_{L_0}a_{I_1}  
			\\&\quad\cdot G_{L_0M_0I_1}^{X_0} G_{I_1M_1L_0}^{X_1}\cdot  u_{[M^\bullet]}.
		\end{align*}
		
		In order to prove the above two expressions are equal to each other, we only need to relate variables in RHS(\ref{dd}) such that $[L_0]=[I_0]$, and $[I_1]=[Z_0]$.
		
		The proof is completed.
	\end{proof}

	
	In \cite{Cr10}, Cramer proved that the reduced Drinfeld double of the extended Hall algebra of a hereditary category is invariant under derived equivalences. We restate the conclusion as follows.
	
	\begin{proposition}
		[\text{\cite[Proposition 5]{Cr10}}]
		\label{cr}
		Let $F:\cd^b(\ca)\rightarrow \cd^b(\mathcal{B})$ be a derived equivalence. Then the following assignment extends to an isomorphism between Drinfeld double $\cd\th(\ca)$ and $\cd\th(\mathcal{B})$:
		\begin{align*}
			\phi:\cd\th(\ca)&\longrightarrow \cd\th(\mathcal{B})\\
			1\otimes [M]&\mapsto \begin{cases}
				1\otimes [N], \text{ for $n$ even},\\
				[N]\otimes 1,\text{ for $n$ odd};
			\end{cases}\\
			[M]\otimes 1&\mapsto \begin{cases}
				[N]\otimes 1, \text{ for $n$ even},\\
				1\otimes [N],\text{ for $n$ odd};
			\end{cases}\\
		\end{align*}
		where $M\in \ca \cap F^{-1}(\mathcal{B}[n])$ and $N=F(M)[-n]\in F(\ca[-n]) \cap \mathcal{B}$.
	\end{proposition}
	
	Let $F:\cd^b(\ca)\rightarrow \cd^b(\mathcal{B})$ be a derived equivalence. Then $F$ induces an equivalence $\overline{F}:\cR(\ca)\rightarrow \cR(\mathcal{B})$. 
	Combining Theorem \ref{iso} and Proposition \ref{cr}, we get the following corollary immediately.
	
	\begin{corollary}
		\label{cor:derinvar}
		For hereditary abelian $\bfk$-categories $\ca$ and $\mathcal{B}$ with Euler form skew symmetric, if there exists a derived equivalence $F:\cd^b(\ca)\rightarrow \cd^b(\mathcal{B})$, then the following assignment extends to an isomorphism between the twisted derived Hall algebras $\widetilde{\ch}(\cR(\ca))$ and $\widetilde{\ch}(\cR(\mathcal{B}))$:
		\begin{align*}
			\Phi:\widetilde{\ch}(\cR(\mathcal{A}))&\longrightarrow \widetilde{\ch}(\cR(\mathcal{B}))\\
			u_{[M]}&\mapsto u_{[\overline{F}M]},
			\\
			u_{[M[1]]}&\mapsto u_{[\overline{F}M[1]]},
		\end{align*}
		for $M\in \ca \cap F^{-1}(\mathcal{B}[n])$. 
	\end{corollary}

	\begin{proof}
		It follows by setting $\Phi=\Psi_\mathcal{B}\circ \phi \circ\Psi^{-1}_\mathcal{A}$.
	\end{proof}
	
	For any two hereditary abelian categories $\ca$ and $\cb$,  we conjecture that their root categories are triangulated equivalent if and only if they are derived equivalent. 
	By Corollary \ref{cor:derinvar}, we know that the derived Hall algebras $\widetilde{\ch}(\cR(\ca))$ are invariant up to derived equivalence. So it is reasonable to call them derived Hall algebras of root categories.

	It is interesting to describe the action of the isomorphism $\Phi:\widetilde{\ch}(\cR(\mathcal{A}))\rightarrow \widetilde{\ch}(\cR(\mathcal{B}))$ on the basis $u_{[M]\oplus [N[1]]}$.

	\section{Derived Hall algebra of the Jordan quiver}
	\label{sec:Jordan}
	
	In this section, we denote by $\QJ$ the Jordan quiver,  i.e., the quiver with a single vertex $1$ and a single loop $\alpha:1\rightarrow 1$.
	Let $\ca={\mathrm {rep}}^\mathrm{nil}_\bfk(\QJ)$ be the category formed by finite-dimensional nilpotent $\bfk$-linear representations. It is well-known that the Euler form of $\ca$ in this case is trivial.
	
	We consider the derived Hall algebra of $\cR({\mathrm {rep}}^\mathrm{nil}_\bfk(\QJ))$ in the following. First, we simplify its structure constants.

	For any $M^\bullet$, $L^\bullet$, $Z^\bullet$ in $\cR(\ca)$, set
	\begin{align*}
		\Hom(Z^\bullet,M^\bullet)_{L^\bullet}&=
		\{l:Z^\bullet\rightarrow M^\bullet\mid \cone l\cong L^\bullet\},
		\\
		\cw(Z^\bullet,L^\bullet;M^\bullet)&=\{(l,m,n)\in\Hom(Z^\bullet,M^\bullet)\times \Hom(M^\bullet,L^\bullet)\times\Hom(L^\bullet,Z^\bullet[1])\mid 
		\\\notag
		&\qquad Z^\bullet \xrightarrow{l}M^\bullet\xrightarrow{m} L^\bullet\xrightarrow{n}Z^\bullet[1] \text{ is a triangle} \}.
	\end{align*}
	Let $  \cv(Z^\bullet,L^\bullet;M^\bullet)$ be the orbit space under the natural action of $\Aut(Z^\bullet)\times \Aut(L^\bullet)$. 
	In this case, we simplify formulas (\ref{3.9})-(\ref{3.11}) as follows:
	\begin{align}
		F_{L^\bullet Z^\bullet}^{M^\bullet}:&=\dfrac{1}{\{Z^\bullet[1],L^\bullet\}}\sum_{\alpha\in	\cv(Z^\bullet,L^\bullet;M^\bullet)}\dfrac{|\End(L^\bullet_1(\alpha))|}{|\Aut(L^\bullet_1(\alpha))|}
		\\&=
		{|\Hom(M^\bullet,L^\bullet)_{Z^\bullet[1]}|}\cdot \dfrac{\{L^\bullet[1],L^\bullet\}   }{|\Aut(L^\bullet)|\{M^\bullet[1],L^\bullet\}}
		\\
		&={|\Hom(Z^\bullet,M^\bullet)_{L^\bullet}|}\cdot \dfrac{\{Z^\bullet[1],Z^\bullet\}   }{|\Aut(Z^\bullet)|\{Z^\bullet[1],M^\bullet\}}.
	\end{align}
	The corresponding Drinfeld dual Hall number (see formula (\ref{3.20})) can also be simplified as
	\begin{align*}
		\dfrac{|\Hom(L^\bullet,Z^\bullet[1])_{M^\bullet[1]}|}{|\Hom(H^0(L^\bullet),H^0(Z^\bullet))||\Hom(H^1(L^\bullet),H^1(Z^\bullet))|},
	\end{align*}
	which can be viewed as the twisted derived Hall number as well.

	\begin{proposition}
		\label{prop:commDHA}
		The twisted derived Hall algebra of $\ca$ is commutative.
	\end{proposition}

	\begin{proof}
		By Lemma \ref{basis}, The twisted derived Hall algebra has a basis $\{u_{[A]}*u_{[B[1]]}\mid A,B\in\ca\}$. So we only need to prove the following three equations for $L,Z\in\ca$:
		\begin{align}
			u_{[L]}*u_{[Z]}&=u_{[Z]}*u_{[L]};\\
			u_{[L[1]]}*u_{[Z[1]]}&=u_{[Z[1]]}*u_{[L[1]]};\\
			u_{[L]}*u_{[Z[1]]}&=u_{[Z[1]]}*u_{[L]}.
		\end{align}
		
		Note that 
		\begin{align*}
			u_{[L]}*u_{[Z]}&=\sum_{[M]}
			\dfrac{|\Hom(L,Z[1])_{M[1]}|}{|\Hom(L,Z)|}\cdot u_{[M]}
			\\\notag&=\sum_{[M]}
			\dfrac{|\Ext^1_\ca(L,Z)_{M}|}{|\Hom(L,Z)|}\cdot u_{[M]}
			\\\notag
			&=\sum_{[M]}G_{LZ}^M \frac{|\Aut(Z)| |\Aut(L)|}{|\Aut(M)|}
			\cdot u_{[M]},
		\end{align*}
		and similarly\[
		u_{[Z]}*u_{[L]}=\sum_{[M]}G_{ZL}^M \frac{|\Aut(Z)| |\Aut(L)|}{|\Aut(M)|}
		\cdot u_{[M]}.
		\]
		It is well-known that $G_{LZ}^M=G_{ZL}^M$; see e.g. \cite[Chaper II, \S2]{Mac95}. Hence we have $u_{[L]}*u_{[Z]}=u_{[Z]}*u_{[L]}$. By a similar argument, we also have $u_{[L[1]]}*u_{[Z[1]]}=u_{[Z[1]]}*u_{[L[1]]}$.
		On the other hand, 
		\begin{align*}
			u_{[L]}*u_{[Z[1]]}&=\sum_{[M^\bullet]} |\Hom(L,Z)_{M^\bullet[1]}|\cdot u_{[M^\bullet]}
			\\
			&=\sum_{[M^\bullet]} |\{l:L\rightarrow Z\mid \ker l\cong M_0,\coker l\cong M_1\}|\cdot u_{[M^\bullet]}
			\\
			&=
			\sum_{[M^\bullet]}\sum_{[I]}G_{IM_0}^{L}G_{M_1I}^{Z}\cdot a_I\cdot u_{[M^\bullet]}
		\end{align*}
		for $M^\bullet\cong M_0\oplus M_1[1]$ with $M_0,M_1\in\ca$. Similarly
		\begin{align*}
			u_{[Z[1]]}*u_{[L]}&=\sum_{[M^\bullet]} |\Hom(Z[1],L[1])_{M^\bullet[1]}|\cdot u_{[M^\bullet]}
			\\
			&=\sum_{[M^\bullet]}\sum_{[I]}G_{M_0I}^{L}G_{IM_1}^{Z}\cdot a_I\cdot u_{[M^\bullet]}.
		\end{align*}
		Due to $G_{IM_0}^L=G_{M_0I}^L$ and $G_{M_1I}^Z=G_{IM_1}^Z$, we get $ u_{[L]}*u_{[Z[1]]}=u_{[Z[1]]}*u_{[L]}$.
	\end{proof}

	Let $\th(\ca)\otimes\th(\ca)$ be the tensor product of two copies of $\th(\ca)$, equipped with the standard multiplication, i.e.,
	{$(a\otimes u)(b\otimes w)=(a b)\otimes (u w)$}. 
	
	\begin{corollary}
		\label{cor:iso1}
		There is an isomorphism of algebras between $\th(\ca)\otimes\th(\ca)$ and $\cd\th(\ca)$.
	\end{corollary}

	\begin{proof}
		There is a natural isomorphism of vector spaces
		\begin{align}
			\varphi:\widetilde{\ch}(\ca)\otimes \widetilde{\ch}(\ca)&\longrightarrow \cd\th(\ca)
			\\
			\notag
			a\otimes b&\mapsto a\otimes b.
		\end{align}
		For any $a,b,c,d\in \widetilde{\ch}(\ca)$, we have 
		\begin{align*}
			\varphi\big((a\otimes b)(c\otimes d)\big)&=\varphi(ac\otimes bd)
			\\
			&=ac\otimes bd
			\\
			&\stackrel{(D2)}{=}(ac\otimes 1)*(1\otimes bd)
			\\
			&\stackrel{(D1)}{=}(a\otimes 1)*(c\otimes 1)*(1\otimes b)*(1\otimes d)
			\\
			&=\big((a\otimes 1)*(1\otimes b)\big)*\big((c\otimes 1)*(1\otimes d)\big)
			\\
			&\stackrel{(D2)}{=}(a\otimes b)*(c\otimes d)
			\\
			&=\varphi(a\otimes b)*\varphi(c\otimes d).
		\end{align*}
		Here the fifth equality follows from Proposition \ref{prop:commDHA}.
		
		Therefore, $\varphi$ is an isomorphism of algebras.
	\end{proof}
	
	Let $\mathfrak{S}_n$ be the symmetric group in $n$ variables for $n\geq1$. 
	Let $\Lambda_n=\mathbb{Q}[x_1,\cdots,x_n]^{\mathfrak{S}_n}$ be the ring of symmetric polynomials in $n$ variables for $n\ge1$. These form a projective system via the maps $\Lambda_{n+1}\rightarrow \Lambda_n$ obtained by setting the last variable $x_{n+1}$ to zero. The projective limit (in the category of graded rings) $\Lambda=\mathop{\lim}\limits_{\longleftarrow}\Lambda_n$ can thus be considered as the ring $\mathbb{Q}[x_1,x_2,\cdots]^{\mathfrak{S}_\infty}$ of
	symmetric polynomials in infinitely many variables.
	It is well-known
	that $\Lambda$ is a polynomial ring $\mathbb{Q}[e_r\mid r\ge1]$, where $e_r$ denotes the $r$-th elementary symmetric function.
	
	The ring $\Lambda$ is also equipped with a canonical coproduct: for $n\ge1$ consider the map $\Delta_n:\Lambda_{2n}\rightarrow \Lambda_n\otimes \Lambda_n$ induced by the embedding
	\begin{align*}
		\mathbb{Q}[x_1,\cdots,x_{2n}]^{\mathfrak{S}_{2n}}&\hookrightarrow \mathbb{Q}[x_1,\cdots,x_{2n}]^{\mathfrak{S}_{n}\times {\mathfrak{S}_{n}}}
		\\&=\mathbb{Q}[x_1,\cdots,x_n]^{\mathfrak{S}_{n}}\otimes \mathbb{Q}[x_{n+1},\cdots,x_{2n}]^{\mathfrak{S}_{n}}
	\end{align*}
	where in the second term the first copy of $\mathfrak{S}_n$ permutes together the variables $x_1,\cdots,x_n$ while the second copy of $\mathfrak{S}_n$ permutes together the variables $x_{n+1},\cdots,x_{2n}$. In the projective limit, the maps $\Delta_n$ give rise to a coproduct $\Delta:\Lambda\rightarrow \Lambda\otimes \Lambda$.
	
	Let $t$ be an indeterminate, and let 
	\begin{align}
		\Lambda_t:=\Lambda\otimes \Q(t).
	\end{align}
	
	It is well known that $\ca={\mathrm {rep}}^\mathrm{nil}_\bfk(\QJ)$ is a uniserial category. 
	Let $S$ be the simple object of $\ca={\mathrm {rep}}^\mathrm{nil}_\bfk(\QJ)$. Then any indecomposable object of $\rep_\bfk^{\rm nil}(\QJ)$  (up to isomorphism) is of the form $S^{(n)}$ of length $n\geq1$.
	Thus the set of isomorphism classes $\Iso(\rep_\bfk^{\rm nil}(\QJ))$ is canonically isomorphic to the set $\cp$ of all partitions, via the assignment
	\begin{align}
		\la =(\la_1,\la_2,\cdots,\la_r)\mapsto S^{(\la)}= S^{(\la_1)}\oplus \cdots \oplus S^{(\la_r)}.
	\end{align}

	\begin{theorem}
		[\text{\cite{Mac95}}] 
		\label{mac}
		There is an isomorphism of bialgebras $\psi_{(q)}:\widetilde{\ch}(\ca)\rightarrow \Lambda_{q^{-1}}$ which maps $[S^{\oplus r}]\mapsto  |\Aut(S^{\oplus r})|\cdot q^{-\frac{r(r-1)}{2}}e_r$, where $S$ is the simple object of $\ca$. Furthermore, $\psi_{(q)}$ transports the nondegenerate Hopf pairing $(,)$ of $\widetilde{\ch}(\ca)$ to the Hall-Littlewood scalar product of $\Lambda_{q^{-1}}$, which is the
		scalar product uniquely determined by the conditions $$\{x,yz\}=\{\Delta(x),y\otimes z\},\qquad
		\{p_r,p_s\}=\delta_{r,s}\dfrac{r}{q^r-1},$$
		where $p_r=\sum_ix^r_i$
		stands for the power sum symmetric function.
	\end{theorem}
	
	By Theorem \ref{mac}, we can define the Drinfeld double of $\Lambda_t$, which is denoted by $\cd\Lambda_t$,  and called the ring of Laurent symmetric functions; see \cite{BS12}. Then the following corollary follows from Corollary \ref{cor:iso1} and Theorem \ref{mac}.
	\begin{corollary}[\text{cf. \cite{Zel81,BS12}}]
		\label{cor5.4}
		There is an isomorphism of algebras from the tensor product algebra $\Lambda_t\otimes\Lambda_t$ to the Drinfeld double $\cd\Lambda_t$ of $\Lambda_t$.
	\end{corollary}

	Combining Theorem \ref{mac}, Corollary \ref{cor5.4} and Theorem \ref{iso}, we have the following corollary.

	\begin{corollary}
		There is an isomorphism of algebras
		\begin{align*}
			\Theta_{(q)}:\widetilde{\ch}(\cR(\ca))&\longrightarrow \cd \Lambda_{q^{-1}},
			\\ u_{[S^{\oplus r}]}&\mapsto  |\Aut(S^{\oplus r})|\cdot q^{-\frac{r(r-1)}{2}} e_r\otimes 1;
			\\ u_{[S^{\oplus r}[1]]}&\mapsto  |\Aut(S^{\oplus r})|\cdot q^{-\frac{r(r-1)}{2}} 1\otimes e_r.
		\end{align*}
	\end{corollary}

	\begin{proof}
		It follows by setting $\Theta_{(q)}=(\psi_{(q)}\otimes \psi_{(q)}) \circ\Psi^{-1}_\mathcal{A}$.
	\end{proof}

	For any partition $\lambda$, let $P_\lambda(x;t)$ be the Hall-Littlewood symmetric function. Then the isomorphism $\psi_{(q)}:\widetilde{\ch}(\ca)\rightarrow \Lambda_{q^{-1}}$ maps $[S^{(\lambda)}]$ to $q^{-\sum\limits_i (i-1)\lambda_i}|\Aut(S^{(\lambda)})|P_\lambda(x;q^{-1})$. In a sequel, we shall introduce Hall-Littlewood Laurent symmetric functions  $P_{\lambda,\mu}(x;t)$ in $\cd\Lambda_t$ for any bipartitions $(\lambda,\mu)$, which form a basis of $\cd\Lambda_t$, and correspond to $u_{[S^{(\lambda)}\oplus S^{(\mu)}[1]]}$ 
	in $\widetilde{\ch}(\cR(\ca))$ via the isomorphism $\Theta_{(q)}: \widetilde{\ch}(\cR(\ca))\rightarrow \cd\Lambda_{q^{-1}}$. In particular, $P_{\lambda,\emptyset}(x;t)=P_\lambda(x;t)=P_{\emptyset, \lambda}(x;t)$ for any partition $\lambda$. These Hall-Littlewood Laurent symmetric functions shall be defined via Giambelli type polynomials, and the (horizontal and vertical) Pieri rules shall be obtained.

	\section{Derived Hall algebras of elliptic curves} 
	\label{sec:elliptic}

	In this section, we shall study the derived Hall algebras of elliptic curves and their connections to double affine Hecke algebras (DAHAs) following \cite{BS12}.
	
	Throughout this section, $\X$ denotes a smooth elliptic curve
	defined over $\bfk$, that is, $\X$ is a smooth projective curve of genus one. 
	We denote by $\ca={\rm coh}(\X)$ its category of coherent sheaves. For a closed point $x$ of $\X$, denote by $\ct_x$ the subcategory of torsion sheaves supported at $x$. There is a unique simple object $\co_x$ in $\ct_x$, and $\ct_x$ is equivalent to the category of finite length modules over the local ring $R_x$ of the point $x$. More explicitly, $\ct_x$ is equivalent to $\rep^{\rm nil}_{\bfk_x}(\QJ)$, where $\bfk_x$ is the residue field $\co_x$, and $\bfk_x$ is a finite field extension of $\bfk$ with $[\bfk_x:\bfk]=\deg(x)$.  
	
	Denote by $\cc_\infty$ the category of torsion sheaves, which is equivalent to the coproduct category $\prod_x\ct_x$, where $x$ runs through the set of closed points of $\X$.  
	
	It is well-known that the Euler form of $\ca={\rm coh}(\X)$ is skew symmetric, that is, the symmetric Euler form is trivial. The slope of a sheaf $\cf$ is $\mu(\cf)=\deg(\cf)/\rank(\cf)$, and 
	a sheaf $\cf$ is semi-stable (resp. stable) if for any subsheaf $\cg\subset \cf$ we have $\mu(\cg)\leq \mu(\cf)$ (resp. $\mu(\cg)<\mu(\cf)$). The full subcategory $\cc_\mu$ consisting of all semi-stable sheaves of a fixed slope $\mu\in \Q\cup\{\infty\}$ is abelian, artinian and closed under extensions. From \cite{A57}, there are canonical exact equivalences of abelian categories 
	\begin{align}
		\epsilon_{\nu,\mu}:\cc_\mu\stackrel{\sim}{\longrightarrow}\cc_\nu.
	\end{align}

	Let $\widetilde{\ch}(\X)$ be the twisted Hall algebra of $\ca={\rm coh}(\X)$. For any closed point of $\X$, let $\th_x$ be the twisted Hall algebra of $\ct_x$, which is a subalgebra of $\th(\X)$. By Theorem \ref{mac}, there is an isomorphism $\psi_x:\th_x \rightarrow \Lambda_{q^{-\deg(x)}}$.
	For any $r\geq1$ we define $T_{r,x}^{(\infty)}\in\th_x$ to be
	\begin{align}
		\frac{T_{r,x}^{(\infty)}}{[r]_\sqq}=\begin{cases} 
			0& \text{ if } \deg(x)\nmid r,
			\\
			\frac{\deg(x)}{r}\psi_x^{-1}(p_{\frac{r}{\deg(x)}}) & \text{ if }\deg(x)\mid r.
		\end{cases}
	\end{align}
	Define $T_r^{(\infty)}=\sum_xT_{r,x}^{(\infty)}$. 
	For any $\mu\in \Q$ we put $T_r^{(\mu)}=\epsilon_{\mu,\infty}(T_r^{(\infty)})$. 

	Let $\cR_\X$ be the root category of $\coh(\X)$, and $\widetilde{\ch}(\cR_\X)$ be the twisted derived Hall algebra of $\coh(\X)$. By Lemma \ref{embed}, we have two embeddings $\Psi^\pm:\th(\X)\rightarrow\widetilde{\ch}(\cR_\X)$. 
	Let $\U_\X\subset \widetilde{\ch}(\cR_\X)$ be the $\Q(\sqq)$-subalgebras   generated by all elements $(T_r^{(\mu)})^\pm:=\Psi^\pm(T_r^{(\mu)})$ for $r\ge 1$ and $\mu\in\Q\cup\{\infty\}$.
	
	If $\mu=\dfrac{l}{n}$ with $n\geq1$ and $l,n$ relatively prime, we put 
	\begin{align}
		T_{(\pm rn,\pm rl)}=(T_r^{(\mu)})^\pm,
	\end{align}
	and
	$T_{(0,\pm r)}=(T_r^{(\infty)})^\pm$ and $T_{(0,0)}=1$. Then $\U_\X$ is the subalgebra of $\widetilde{\ch}(\cR_\X)$ generated by $T_{(m,n)}$ for $(m,n)\in\Z^2\setminus \{(0,0)\}$. 
	
	For any $\bx=(m,n)\in\Z^2\setminus \{(0,0)\}$, denote by $\deg(\bx)=\gcd(m,n)$.
	
	For formal parameters $\sigma$ and $\overline{\sigma}$, let $\mathbf{R}=\C[\sigma^{\pm1/2},\overline{\sigma}^{\pm1/2}]$,
	Burban and Schiffmann \cite[\S5.2,\S6.2]{BS12} also introduced a $\C$-algebra $\ce_\mathbf{R}$ generated by $\{\tt_\bx\mid \bx\in\Z^2\setminus \{(0,0)\}\}$ subject to some relations, and we do not recall its definition here.
	
	Let $\sharp \X(\F_{q^r})$ be the number of rational points of $\X$ over $\F_{q^r}$. Then there exist conjugate algebraic numbers $\sigma,\overline{\sigma}$, satisfying $\sigma\overline{\sigma}=q$, such that
	\begin{align}
		\label{eq:number}
		\sharp\X(\F_{q^r})=q^r+1-(\sigma^r+\overline{\sigma}^r)
	\end{align}
	for any $r\geq1$.

	\begin{proposition}
		If $\sigma,\overline{\sigma}$ are algebraic numbers satisfying \eqref{eq:number}, and we denote by $\ce_{\sigma,\overline{\sigma}}$ the $\C$-algebra obtained from $\ce_\mathbf{R}$ naturally, then the assignment $\Omega:\tt_\bx\mapsto [\deg(\bx)]_\sqq T_\bx$ for $\bx\in\Z^2\setminus\{(0,0)\}$ extends to an isomorphism $\Omega:\ce_{\sigma,\overline{\sigma}}\stackrel{\simeq}{\longrightarrow} \U_\X\otimes \C$.
	\end{proposition}
	
	\begin{proof}
		It follows from \cite[Theorem 5.4]{BS12} and Theorem \ref{iso}.
	\end{proof}

	We can view $\ce_\mathbf{R}$  as a generic version of $\U_\X$.

	Let \begin{align}
		\bM=\C[x_1^{\pm1},x_2^{\pm1},\dots,y_1^{\pm1},y_2^{\pm1},\dots]^{\mathfrak{S}_\infty}
	\end{align}
	be the ring of diagonal symmetric polynomials,
	with $\mathfrak{S}_\infty$ acting simultaneously on the variables $x_i$ and $y_i$. It is proved in \cite{BS12} that there is an isomorphism 
	\begin{align}
		(\ce_\mathbf{R})|_{\sigma=\overline{\sigma}=1}\cong \bM.
	\end{align}
	
	The spherical DAHA $\mathbf{S}{\textbf{\" H}}_n$ of type $\mathfrak{gl}(n)$ is introduced in \cite{Ch05}, which is a two-parameter deformation of the ring
	\begin{align}
		\bM_n=\C[x_1^{\pm1},\dots,x_n^{\pm1},y_1^{\pm1},\dots,y_n^{\pm1}]^{\mathfrak{S}_n}.
	\end{align}
	It is proved in \cite{SV1} that there are surjective homomorphisms $\ce_\mathbf{R} \rightarrow  \mathbf{S}{\textbf{\" H}}_n$ for any positive
	integer $n$, so that $\ce_\mathbf{R}$ can be thought of as the “stable limit” $\mathbf{S}{\textbf{\" H}}_\infty$ of the type
	A spherical DAHA. 
	
	Compared with the Drinfeld double of Hall algebra used in \cite{BS12}, the twisted derived Hall algebras $\widetilde{\ch}(\cR_\X)$ and $\U_\X$ are intrinsic and categorical.  
	In a sequel, we shall study the algebra $\U_\X$ defined here in detail, and make its connection to the DAHA more clearly.
	
	\appendix
	
	\section{}

	In this section, we shall give two lemmas for arbitrary abelian categories $\cb$. These lemmas are used to prove Proposition \ref{prop1}.
	
	\begin{lemma}
		\label{lem1}
		Given a commutative diagram of short exact sequences in an abelian category $\cb$ as follows
		\[\begin{tikzpicture}
			\node (-2) at (-2,0) {$0$};
			\node (2) at (0,0) {$V_1'$};
			\node (8) at (7,0) {$0$};
			\node (81) at (7,-2) {$0$};
			
			\draw[->] (-2) --node[above ]{} (2);
			\node (3) at (2.5,0) {$V_2$};
			\draw[->] (2) --node[above ]{$f'$} (3);
			\node (4) at (3+1+1,0) {$V_3$};
			\draw[->] (3) --node[above ]{$g$} (4);
			
			\draw[-] (2.5,-.3) --node[above ]{} (2.5,-1.7);
			\draw[-] (2.6,-.3) --node[above ]{} (2.6,-1.7);
			
			\node (-21) at (-2,-2) {$0$};
			\node (21) at (0,-2) {$V_1$};
			\draw[->] (-21) --node[above ]{} (21);
			\node (31) at (2.5,-2) {$V_2$};
			\draw[->] (21) --node[above ]{$f$} (31);
			\node (41) at (3+1+1,-2) {$V_3'$};
			\draw[->] (31) --node[above ]{$g'$} (41);
			\draw[->] (4) --node[above ]{} (8);
			
			\draw[->] (41) --node[above ]{} (81);
			\draw[->] (2) --node[right ]{$h'$} (21);
			
			\draw[->] (4) --node[right ]{$h$} (41);
		\end{tikzpicture}\]
		we have
		\begin{align*}
			\widehat{V_1}+\widehat{V_3}=\widehat{V_2}+\widehat{\Im (g\circ f)}
		\end{align*}
		in $K_0(\cb)$. 
		Moreover, if $\cb$ is the category of finite-dimensional $\bfk$-linear spaces, then
		\begin{align*}
			|V_1|\cdot|V_3|=|V_2|\cdot |\Im(g\circ f)|.
		\end{align*}
	\end{lemma}
	
	\begin{proof}
		By assumption $f'=f\circ h'$ and $f'$ is injective, we have $h'$ is injective. Then we obtain the following commutative diagram:
		\[\begin{tikzpicture}
			\node (-2) at (-2,0) {$0$};
			\node (2) at (0,0) {$V_1'$};
			\node (8) at (7,0) {$0$};
			\node (81) at (7,-2) {$0$};
			
			\draw[->] (-2) --node[above ]{} (2);
			\node (3) at (2.5,0) {$V_1$};
			\draw[->] (2) --node[above ]{$h'$} (3);
			\node (4) at (3+1+1,0) {$\coker h'$};
			\draw[->] (3) --node[above ]{$\pi$} (4);
			
			\draw[-] (0,-.3) --node[above ]{} (0,-1.7);
			\draw[-] (0.1,-.3) --node[above ]{} (0.1,-1.7);
			
			\node (-21) at (-2,-2) {$0$};
			\node (21) at (0,-2) {$V_1'$};
			\draw[->] (-21) --node[above ]{} (21);
			\node (31) at (2.5,-2) {$V_2$};
			\draw[->] (21) --node[above ]{$f'$} (31);
			\node (41) at (3+1+1,-2) {$V_3$};
			\draw[->] (31) --node[above ]{$g$} (41);
			\draw[->] (4) --node[above ]{} (8);
			
			\draw[->] (41) --node[right ]{} (81);
			\draw[->] (3) --node[right ]{$f$} (31);
			
			\draw[->] (4) --node[right ]{$i$} (41);
		\end{tikzpicture}\]
		which yields an exact sequence
		$$0\longrightarrow V_1\longrightarrow V_2\oplus \coker h'\longrightarrow V_3\longrightarrow 0.$$
		Note that $\pi$ is surjective and $i$ is injective with $g\circ f=i\circ \pi$. Hence $\Im(g\circ f)\cong \coker h'$. Then the result follows.
	\end{proof}

	\begin{proposition}
		\label{cor1}
		Given morphisms $f:V_1\rightarrow V_2$ and $g:V_2\rightarrow V_3$ in an abelian category $\cb$ satisfying $p\circ i=0$, where  $(C,p)$ is the cokernel of $f$ and $(K,i)$ is the kernel of $g$, we have
		\begin{align*}
			\widehat{\Im(f)}+\widehat{\Im(g)}=\widehat{V_2}+\widehat{\Im(g\circ f)}.
		\end{align*}
		Moreover, if $\cb$ is the category of finite-dimensional $\bfk$-linear spaces, then
		\begin{align*}
			|\Im  (f)|\cdot  |\Im  (g)|=|V_2|\cdot |\Im  (g\circ f)|.
		\end{align*}
	\end{proposition}
	
	\begin{proof}
		Assume $f=i_f\circ p_f$ with $i_f:\Im (f)\rightarrow V_2$ a monomorphism and $p_f:V_1\rightarrow \Im (f)$ an epimorphism, $g=i_g\circ p_g$ with $i_g:\Im (g)\rightarrow V_3$ a monomorphism and $p_g:V_2\rightarrow \Im (g)$ an epimorphism. Then we have the following commutative diagram
		\[\xymatrix{V_1 \ar[dr]_{f}\ar[r]^-{p_f} & \Im f \ar[r]^-{p'} \ar[d]^{i_f} & I\ar[d]^-{i'} \\
			&V_2\ar[dr]_{g}\ar[r]^-{p_g} &\Im g \ar[d]^-{i_g} \\
			&& V_3}\]
		where $I$ is the image of $p_g\circ i_f$, and $p',i'$ are the canonical morphisms. 
		
		Note that $C$ is the cokernel of $i_f$, $K$ is the kernel of $p_g$ and $I$ is also the image of $g\circ f$. We get the following commutative diagram
		\[\xymatrix{ 0\ar[r] & K\ar[r]^i\ar[d] & V_2\ar[r]^{p_g} \ar@{=}[d] & \Im g \ar[r]\ar[d] &0
			\\
			0\ar[r]& \Im f \ar[r]^{i_f} &V_2\ar[r]^{p} &C\ar[r] &0 }\]
		The desired results follow from  Lemma \ref{lem1}.
	\end{proof}


\end{document}